\documentclass[11pt]{article}

\usepackage{amsmath, amsfonts, amsthm, verbatim, amssymb}

\newcommand{\dd}{\mathrm{d}}
\newcommand{\R}{\mathbb{R}}
\newcommand{\N}{\mathbb{N}}
\newcommand{\E}{\mathbb{E}}

\newtheoremstyle{neu-theorem}
{11pt}      
{11pt}      
{\itshape}                  
{}          
{\bfseries} 
{}          
{1em}  
{\textbf{\thmname{#1}\thmnumber{ #2}\thmnote{ (#3)}}}          

\theoremstyle{neu-theorem}

\newtheorem{theorem}{Theorem}[section]
\newtheorem{corollary}[theorem]{Corollary}
\newtheorem{lemma}[theorem]{Lemma}
\newtheorem{proposition}[theorem]{Proposition}

\newtheoremstyle{neu}
{11pt}      
{11pt}      
{}                  
{}          
{\bfseries} 
{}          
{1em}  
{\textbf{\thmname{#1}\thmnumber{ #2}\thmnote{ (#3)}}}          

\theoremstyle{neu}
\newtheorem{definition}[theorem]{Definition}
\newtheorem{example}[theorem]{Example}
\newtheorem{remark}[theorem]{Remark}
\newtheorem{assumption}{Assumption}

\usepackage{latexsym}
\usepackage{amssymb}
\usepackage{amsmath}
\usepackage{amsthm}
\usepackage{mathrsfs}
\usepackage{url}
\usepackage{longtable}
\usepackage{color}
\usepackage{dsfont}
\usepackage[T1]{fontenc}
\usepackage{lmodern}
\usepackage[linktocpage]{hyperref}

\usepackage{graphicx}
\usepackage[latin1]{inputenc}
\usepackage{amsopn}

\newcommand{\bfi}{\begin{fig}}
	\newcommand{\efi}{\end{fig}}
\newcommand{\btab}{\begin{tab}}
	\newcommand{\etab}{\end{tab}}
\newcommand{\barr}{\begin{array}}
	\newcommand{\earr}{\end{array}}
\newcommand{\beq}{\begin{equation}}
	\newcommand{\eeq}{\end{equation}}
\newcommand{\bdis}{\begin{displaymath}}
	\newcommand{\edis}{\end{displaymath}\noindent}

\newcommand{\bbn}{\mathbb{N}}
\newcommand{\bbz}{\mathbb{Z}}
\newcommand{\bbr}{\mathbb{R}}

\newcommand{\bbe}{\mathbb{E}}

\newcommand{\bbp}{\mathbb{P}}

\newcommand{\bone}{\mathds 1}

\newcommand{\calf}{{\cal F}}

\newcommand{\calb}{{\cal B}}

\newcommand{\calg}{{\cal G}}

\newcommand{\al}{{\alpha}}
\newcommand{\la}{{\lambda}}
\newcommand{\La}{{\Lambda}}
\newcommand{\eps}{{\epsilon}}

\newcommand{\ga}{{\gamma}}
\newcommand{\Ga}{{\Gamma}}

\newcommand{\si}{{\sigma}}
\newcommand{\om}{{\omega}}
\newcommand{\Om}{{\Omega}}

\newcommand{\ov}{\overline}
\newcommand{\un}{\underline}

\newcommand{\ee}{\mathrm{e}}

\newcommand{\Comp}{\mathrm{c}}

\definecolor{orange}{rgb}{1,0.5,0}

\newcommand{\bthm}{\begin{theorem}}
	\newcommand{\ethm}{\end{theorem}}

\newcommand{\bcor}{\begin{corollary}}
	\newcommand{\ecor}{\end{corollary}}

\newcommand{\blem}{\begin{lemma}}
	\newcommand{\elem}{\end{lemma}}

\newcommand{\bprop}{\begin{proposition}}
	\newcommand{\eprop}{\end{proposition}}

\newcommand{\bdf}{\begin{definition}}
	\newcommand{\edf}{\end{definition}}

\newcommand{\bex}{\begin{example}}
	\newcommand{\eex}{\end{example}}

\newcommand{\brem}{\begin{remark}}
	\newcommand{\erem}{\end{remark}}

\newcommand{\bass}{\begin{assumption}}
	\newcommand{\eass}{\end{assumption}}

\newcommand{\bpr}{\begin{proof}}
	\newcommand{\epr}{\end{proof}}

\newcommand{\benu}{\begin{enumerate}}
	\newcommand{\eenu}{\end{enumerate}}

\newcommand{\bit}{\begin{itemize}}
	\newcommand{\eit}{\end{itemize}}

\newcommand{\bff}{\textbf}

\numberwithin{equation}{section}

\title{Intermittency for the stochastic heat equation with Lévy noise}

\author{
Carsten Chong\thanks{Center for Mathematical Sciences, Technical University of Munich, 
Boltzmannstraße 3, 85748 Garching, Germany, e-mail: carsten.chong@tum.de} \and Péter Kevei\thanks{{MTA-SZTE Analysis and Stochastics 
Research Group, Bolyai Institute, Aradi vértanúk tere 1, 6720 Szeged, Hungary,
e-mail: kevei@math.u-szeged.hu}}
}

\begin{document}
 
 \date{}
\maketitle

\begin{abstract} \noindent
We investigate the moment asymptotics of the solution to the stochastic heat equation 
driven by a $(d+1)$-dimensional Lévy space--time white noise. Unlike the case of Gaussian 
noise, the solution typically has no finite moments of order $1+2/d$ or higher. 
Intermittency of order $p$, that is, the exponential growth of the $p$th moment as time 
tends to infinity, is established in dimension $d=1$ for all values $p\in(1,3)$, and in 
higher dimensions for some $p\in(1,1+2/d)$. The proof relies on a new moment lower bound 
for stochastic integrals against compensated Poisson measures. The behavior of the 
intermittency exponents when $p\to 1+2/d$ further indicates that intermittency in the 
presence of jumps is much stronger than in equations with Gaussian noise. The effect of 
other parameters like the diffusion constant or the noise intensity on intermittency will 
also be analyzed in detail.
\end{abstract}

\vfill

\noindent
\begin{tabbing}
	{\em AMS 2010 Subject Classifications:} \= primary: \,\,\,\,\,\, 60H15, 37H15 \\
	\> secondary: \,\,\, 60G51, 35B40
\end{tabbing}

\vspace{1cm}

\noindent
{\em Keywords:}
comparison principle; intermittency; intermittency fronts; Lévy noise; moment Lyapunov exponents; stochastic heat equation; stochastic PDE

\vspace{0.5cm}

\newpage

\section{Introduction}

We consider the stochastic heat equation on $\bbr^d$ given by
\beq\label{SHE}
\begin{split}
	\partial_t Y(t,x) &= \frac{\kappa}{2}\Delta Y(t,x) + \si(Y(t,x))\dot\La (t,x),\quad 
	(t,x)\in(0,\infty)\times\bbr^d,\\
	Y(0,\cdot)&= f,
\end{split}
\eeq 
where $\kappa\in(0,\infty)$ is the diffusion constant, $\si$ a globally Lipschitz function and $f$ a bounded measurable function on 
$\bbr^d$. The forcing term $\dot \La$ that acts in a multiplicative way on the right-hand 
side of \eqref{SHE} is a \emph{Lévy space--time white noise}, which is the 
distributional derivative of a \emph{Lévy sheet} in $d+1$ parameters. More precisely, we 
assume that $\La$ takes the form
\begin{equation} \label{eq:Lambda}
	\begin{split}
		\Lambda(\dd t, \dd x) 
		& = b \, \dd t\, \dd x + \rho\, W(\dd t, \dd x) + \int_\R z \, 
		(\mu - \nu)(\dd t, \dd x, \dd z),
	\end{split}
\end{equation}
where $b\in\R$ is the mean of $\La$, $\rho \in \R$ is the Gaussian part of $\La$, $W$ is a Gaussian space--time white noise (see \cite{Walsh86}), 
$\mu$ is a Poisson measure on $(0,\infty) \times \R^d \times \R$ with intensity measure
$\nu(\dd t, \dd x, \dd z) = \dd t \,\dd x\, \lambda(\dd z)$, and $\la$ is a Lévy 
measure satisfying
\[ 
\la(\{0\}) =0\quad \text{and} \quad \int_\bbr (1\wedge |z|^2)\,\la(\dd z)<\infty.
\] 
Under the assumption that there exists $p\in[1,1+2/d)$ with
\begin{equation} \label{eq:mlambda}
	m_\lambda( p ) := \left(\int_\R |z|^p\, \lambda(\dd z)\right)^{\frac{1}{p}} <\infty,
\end{equation}
it is shown in \cite{SLB98}  that \eqref{SHE} admits a unique \emph{mild solution} $Y$ satisfying
\beq\label{local-bound} \sup_{(t,x)\in[0,T]\times\R^d} \|Y(t,x)\|_p =  \sup_{(t,x)\in[0,T]\times\R^d} \bbe[|Y(t,x)|^p]^{\frac{1}{p}} < \infty \eeq
for all $T\geq0$. A mild solution to \eqref{SHE} is a predictable process $Y$ satisfying the stochastic Volterra equation 
\begin{equation} \label{eq:heat}
	Y(t,x) = Y_0(t,x) + \int_0^t \int_{\R^d} g(t-s, x-y) \sigma(Y(s,y)) \,\Lambda(\dd s, \dd y),\quad (t,x)\in(0,\infty)\times\bbr^d,
\end{equation}
where 
\beq\label{Y0} 
Y_0(t,x):= \int_{\bbr^d} g(t,x-y)f(y)\,\dd y,\quad (t,x)\in(0,\infty)\times\bbr^d,
\eeq
and 
\begin{equation} \label{eq:heatkernel}
	g(t,x) := g(\kappa; t,x):= \frac{1}{(2\pi\kappa  t)^{d/2}} e^{- \frac{|x|^2}{2\kappa t}}, \quad
	(t,x)\in (0,\infty)\times \R^d,
\end{equation}
is the heat kernel in dimension $d$. As proved in \cite{Chong1}, condition \eqref{eq:mlambda} can be relaxed to include Lévy noises with bad moment properties such as $\al$-stable noises, but in this paper, we will work with \eqref{eq:mlambda} as a standing assumption. 

Our goal is to investigate the behavior of the moments of the solution $Y$ as time tends 
to infinity. In particular, we are interested in conditions under which the solution $Y$ 
to \eqref{eq:heat} exhibits the phenomenon of intermittency. The following definition follows \cite{Carmona94}, Definition~III.1.1, \cite{Conus12}, Equations~(1.6) and 
(1.7), and \cite{Khos}, Definition 7.5.
\bdf\label{def:interm} Let $Y$ be the mild solution to \eqref{eq:heat} and 
$p\in(0,\infty)$. 
\benu
\item $Y$ is said to be \emph{weakly intermittent of order $p$} if
\beq\label{def:interm-2} 0<\un\ga(p)\leq \ov\ga(p)<\infty,\eeq
where the \emph{lower} and \emph{upper moment Lyapunov exponents} $\un\ga(p)$ and $\ov\ga(p)$ are defined as
\begin{align}\label{upperlower} 
\un \ga(p)&:= \liminf_{t\to\infty} \frac{1}{t} \inf_{x\in\bbr^d} \log\bbe[|Y(t,x)|^p]\quad\text{and}\quad
\ov \ga(p):= \limsup_{t\to\infty} \frac{1}{t} \sup_{x\in\bbr^d} \log\bbe[|Y(t,x)|^p]. 
\end{align}
\item $Y$ is said to have a \emph{linear intermittency front of order $p$} if
\beq\label{def:intfronts} 0<\un\la(p)\leq \ov\la(p)<\infty, \eeq
where the \emph{lower} and \emph{upper intermittency fronts} $\un\la(p)$ and $\ov\la(p)$ are defined as
\beq
\begin{split}
		\un\la(p)&:=\sup\left\{ \al>0\colon \limsup_{t\to\infty} \frac{1}{t} \sup_{|x|\geq\al t} 
	\log \bbe[|Y(t,x)|^p]>0 \right\},\\
	\ov\la(p)&:= \inf \left\{ \al > 0 \colon \limsup_{t\to\infty} \frac{1}{t}
	\sup_{|x|\geq \al t} \log \bbe[|Y(t,x)|^p] < 0 \right\},
\end{split}
\eeq
with the convention that $\sup \emptyset :=0$ and $\inf\emptyset:=+\infty$.
\eenu
\edf

For important classes of random fields, the purely moment based notion of weak intermittency in \eqref{def:interm-2} translates into an interesting path property called \emph{physical intermittency}: With high probability, the random field exhibits an extreme mass concentration at large times, in the sense that it almost vanishes on $\bbr^d$ except for exponentially small areas where it develops a whole cascade of exponentially sized peaks. We refer to \cite[Section~2.4]{Bertini95} for a precise statement. 

Similarly, if the initial condition $f$ decays at infinity (in this case we cannot expect to have \eqref{def:interm-2} because of lacking uniformity in the spatial variable), the property \eqref{def:intfronts} would indicate that intermittency peaks, originating from the initial mass around the origin, spread in space at a (quasi-)linear speed.

\subsection*{Review of literature}

The intermittency problem has been investigated by many authors in various situations. For example, \cite{Carmona94} is a classical reference for intermittency in the \emph{parabolic Anderson model (PAM)} on $\bbz^d$, which is the discrete-space analogue of \eqref{SHE} with
\beq\label{PAM} \si(x) = \si_0x, \quad x\in\R, \eeq
for some $\si_0>0$. For the stochastic heat equation, and in particular the continuous PAM driven by a Gaussian space--time white noise, this is analyzed in all its facets in \cite{Bertini95, Chen15, Conus12, Foondun09, Foondun10}, just to name a few. We also refer to \cite{Khos} for a good overview of the subject.

When it comes to stochastic PDEs with non-Gaussian noise, there is much less literature on this topic. Apart from work on the discrete PAM (see \cite{Ahn92, Cranston05} and the references therein), we are only aware of \cite{Balan16} that considers the intermittency problem in continuous space and time. This article investigates the Lévy-driven \emph{stochastic wave equation} in one spatial dimension, and shows that the solution is weakly intermittent of any order $p\geq2$ under natural assumptions. For the proof of the intermittency upper bounds, the authors employ predictable moment estimates for Poisson stochastic integrals, which are surveyed in \cite{MR} in detail. The proof of the lower bound, by contrast, relies on $L^2$-techniques, which are the same as in the Gaussian case treated in \cite{Dalang09} or \cite{Khos}.

\subsection*{Summary of results}

For the stochastic heat equation \eqref{SHE}, however, there is an 
important difference that necessitates the development of new techniques for the 
intermittency analysis. Namely, as soon as $\La$ contains a non-Gaussian part, the 
solution to \eqref{SHE} will typically have finite moments only up to the order 
$(1+2/d)-\epsilon$, even if $\La$ itself has moments of all orders or has bounded jump 
sizes like in the case of a standard Poisson noise, see Theorem~\ref{nomoments}. In particular, 
as soon as we are in dimension $d\geq2$, the solution has \emph{no} finite second moment. This is 
in sharp contrast to the Gaussian case where it is well known that the solution to the 
stochastic heat equation, if it exists, has finite moments of all orders.  
And because, as a consequence of the \emph{comparison principle} in Theorem~\ref{compprin}, we cannot expect in general that the 
solution is weakly intermittent of order $1$, we are forced to consider moments of 
\emph{non-integer} orders in the range $(1,1+2/d)\subseteq (1,2)$. Therefore, 
well-established techniques for estimating \emph{integer} moments of the solution (see 
\cite{Bertini95, Chen15}) do not apply in this setting.

This problem can be remedied by an appropriate use of the 
Burkholder--Davis--Gundy (BDG) inequalities for verifying the intermittency \emph{upper} 
bounds, see Theorem~\ref{upper-bounds}. However, for the corresponding \emph{lower} 
bounds, the moment estimates that are available in the literature (including again the BDG 
inequalities, but also ``predictable'' versions thereof, see e.g.\ \cite{MR}) do not 
combine well with the recursive Volterra structure of \eqref{eq:heat}. So although these 
estimates are sharp, we cannot apply them to produce the desired intermittency lower 
bounds. In order to circumvent this, we use decoupling techniques to establish an -- up to 
our knowledge -- new moment lower bound for Poisson stochastic integrals in 
Lemma~\ref{Poi-ineq}, which we think is of independent interest. With this inequality we 
then prove the weak intermittency of \eqref{SHE} under quite general assumptions. More 
precisely, if $\La$ has mean zero, we show in Theorem~\ref{interm} and 
Theorem~\ref{intfronts} that we have $p$th order intermittency for all $p\in(1,3)$ in 
dimension $1$, and for some $p\in(1,1+2/d)$ in dimensions $d\geq2$. In the latter case, a 
small diffusion constant $\kappa$, or a high noise intensity also leads to intermittency 
of any desired order. Noises with positive or negative mean are treated in 
Theorem~\ref{pos-mean} or Theorem~\ref{interm-neg}, respectively. Moreover, the moment 
estimates in Lemma~\ref{Poi-ineq} also permit us to determine the asymptotics of the 
intermittency exponents as $p\to1+2/d$ or $\kappa\to0$, see Theorem~\ref{asymptotics-1}. 
The results suggest that intermittency in the Lévy case is much more pronounced than 
with Gaussian noise.

Our proofs further indicate that the principal source of intermittency is 
different between the jump and the Gaussian case. In fact, intermittency in the Gaussian 
case is caused by the slow decrease in time of the heat kernel, so peaks in the past are 
remembered for a long time and accumulate to new peaks in the future. By contrast, in the 
Lévy-driven equation, it is the singularity of the heat kernel at the origin that causes 
the high-order intermittent behavior of the solution. So here, for $p$ close to $1+2/d$, 
peaks of order $p$ amplify over short time and hence generate even higher peaks. We refer 
to Remark~\ref{origin} for details.

In the sequel, we will use the letter $C$ to denote a constant whose value may change from line to line and does not depend on anything important in the given context. Sometimes, if we want to stress the dependence of the constant on an important parameter, say $p$, we will write $C_p$. Furthermore, for reasons of brevity, we write $\iint_a^b$ and $\iiint_a^b$ for $\int_a^b\int_{\R^d}$ and $\int_a^b\int_{\R^d}\int_\R$, respectively.



\section{Intermittency upper bounds}\label{Sect2}

We first investigate the upper indices $\ov\ga(p)$ and $\ov\la(p)$, respectively. For a random field $\Phi(t,x)$, indexed by $(t, x)\in (0,\infty)\times \R^d$, and exponents $\beta \in\R$, $c\in[0,\infty)$ 
and $p \in [1,\infty)$, we use the notation 
\begin{equation} \label{eq:norm}
	\|\Phi\|_{p, \beta,c} := \sup_{t \in (0,\infty)} \sup_{x \in \R^d} e^{-\beta t+c|x|} 
	\| \Phi(t,x) \|_p,
\end{equation}
and 
\begin{equation} \label{eq:stochconv}
	( g \circledast \Phi)(t,x) = \int_0^t \int_{\R^d} g(t-s, x-y) \Phi(s,y)\, \Lambda( \dd s, \dd 
	y)
\end{equation}
if $\Phi$ is predictable and the stochastic integral \eqref{eq:stochconv} exists for all $(t,x)\in(0,\infty)\times\R^d$.
The key ingredient for the intermittency upper bounds is the following $L^p$-estimate for 
stochastic convolutions. The Gaussian case has been obtained in 
\cite[Proposition~2.5]{Conus12} and \cite[Proposition~5.2]{Khos}.

\begin{proposition}[Weighted stochastic Young inequality] \label{prop:stochYoung-d}
	Let $d\in\bbn$, $1 \leq p < 1 + 2/d$ and assume that $\rho =0$ if $p<2$.
	For any $c\geq0$ and $\beta>\kappa  c^2 d/2$, we have
	\begin{equation} \label{eq:stochYoung-ineq-d}
		\begin{split}
			\| g \circledast \Phi\|_{p,\beta,c} &\leq 
			C_{\beta,c}(\kappa, p) \, \|\Phi\|_{p,\beta,c}
		\end{split}
	\end{equation}
	with
	\beq\label{Cbetagamma} \begin{split}
	C_{\beta,c}(\kappa,p) &= 
	C_p\Bigg(\frac{2^d |b|}{\beta-\frac12 \kappa c^2 d}
	+  
\frac{2^{\frac{d(3-p)}{2p}}\Gamma(1-\frac{d}{2}(p-1))^\frac{1}{p}m_\la(p)}{p^{\frac{
2+(2-p)d}{2p}}(\pi\kappa)^{\frac{d(p-1)}{2p}} (\beta-\frac12 \kappa c^2 
		d)^{\frac{2-d(p-1)}{2p}}}
	\\
	&\quad+  \frac{m_\la(2)+|\rho|}{(2\kappa(\beta- \frac12 \kappa c^2))^{\frac14}} \bone_{\{d=1,\ p\geq2\}} 
	\Bigg),
	\end{split}\eeq
	where $C_p > 0$ does not depend on $\Lambda$, $\kappa$, $\beta$, $c$ or $d$, 
and it is bounded on $[1+\epsilon, 1 + 2/d)$ for any $\epsilon > 0$.
\end{proposition}

The assumption in Proposition~\ref{prop:stochYoung-d} that $\rho=0$ if $p<2$ means that if $d\geq2$, then necessarily the Gaussian part vanishes because $p<1+2/d\leq 2$. This is reasonable since the stochastic heat equation \eqref{SHE} has no function-valued solution in general if $d\geq2$ and $\rho>0$, see e.g.\ \cite[Section~3.5]{Khos}. Moreover, in dimension $d=1$, we shall only consider the case $p\geq2$ if $\rho>0$. The reason behind is that in the case of Gaussian noise, intermittency of order less than $2$ is open, see the remark after Theorem~\ref{interm}.

\brem The three terms in \eqref{Cbetagamma} illustrate in a nice way the different 
contributions of the noise to the size of $g\circledast \Phi$. The first part comes from 
the deterministic drift  of the noise, the second summand is the $L^p$-contribution 
originating from the jumps, and the third term is the $L^2$-contribution of the jumps and 
the Gaussian part (if $p\geq2$). It is important to notice that a Gaussian noise alone has 
no extra $L^p$-contribution to $C_{\beta,c}(\kappa,p)$ for $p>2$, which reflects the 
equivalence of moments of the normal distribution. Furthermore, as $p \to 1+2/d$, the 
second term explodes for all non-trivial Lévy measures $\la$, no matter how good their 
integrability properties are. This is a first indication that the solution to a 
Lévy-driven stochastic heat equation \eqref{SHE} usually has no finite moments of order 
$1+2/d$ or higher. We confirm this rigorously in Theorem~\ref{nomoments} below.
\erem

With the help of Proposition~\ref{prop:stochYoung-d}, we can extend  the local moment 
bound \eqref{local-bound} obtained in \cite{SLB98} to a global bound. 

\bprop \label{thm:moments}
Assume that 
$f$ satisfies $|f(x)|=O(e^{-c |x|})$ as $|x|\to\infty$ for some $c\geq0$ and that 
$\sigma$ in \eqref{eq:heat} is Lipschitz continuous with
\[
| \sigma(x) -  \sigma(y)| \leq L |x-y|,  \quad x,y \in \R,
\]
for some $L>0$, and also $\si(0)=0$ if $c>0$. Further suppose that $\La$ takes the form \eqref{eq:Lambda} and satisfies \eqref{eq:mlambda} for some $1\leq p < 1+2/d$ as well as $\rho=0$ if $p<2$. Then there 
exists a number $\beta_0>0$ such that the 
stochastic heat equation \eqref{eq:heat} has a unique mild solution $Y$ (up to modifications) 
with $\|Y\|_{p,\beta,c}<\infty$ for all $\beta \geq \beta_0$.
\eprop

We obtain as an immediate consequence upper bounds for the moments of the solution $Y$ to 
the stochastic heat equation \eqref{eq:heat}.
\bthm[Intermittency upper bounds]\label{upper-bounds} Grant the assumptions and notations of 
Proposition~\ref{thm:moments}.
\benu 
\item We have $\ov \ga(p) < \infty$.
\item If $c>0$ and $\si(0)=0$, then
$\ov\la(p)< \infty$.
\eenu
\ethm

\section{Intermittency lower bounds}\label{Sect3}

\subsection{High moments}

One important difference between the stochastic heat equation with jump noise and with Gaussian noise is that the solution $Y$ to \eqref{eq:heat} has no large moments, even in 
dimension $d=1$ and no matter how good the integrability properties of the jumps are. In 
order to understand this, let us consider the situation where $\si\equiv 1$, $f\equiv0$, 
and $\La$ is a standard Poisson random measure, that is, $\la=\delta_1$, $b=1$, and 
$\rho=0$. 
Denoting by $(S_i,Y_i)$ the space--time locations of the jumps of $\La$, we have for 
$(t,x)\in(0,\infty)\times\bbr^d$,
\begin{align*}
	Y(t,x) = \int_0^t \int_{\R^d} g(t-s, x-y)\,  \La(\dd s, \dd y) &= \sum_{i=1}^\infty 
	g(t-S_i,x-Y_i)\bone_{\{S_i<t\}}. 
\end{align*}
If $t>1$, conditionally on the event that at least one point falls into  
$(t-1,t)\times \prod_{i=1}^d (x_i-1,x_i)$, we have
\[
Y(t,x)\geq  g( U, V)= \frac{1}{ (2 \pi\kappa U)^{\frac{d}{2}}} 
e^{-\frac{|V|^2}{2\kappa U}},
\]
where $U,V_1, \ldots, V_d$ are independent and uniformly distributed on $(0,1)$, and $V = 
(V_1, \ldots, V_d)$. Now
\begin{align*}
	\E [g(U,V)^p]  &= \frac{1}{(2 \pi\kappa)^{\frac{pd}{2}}} \int_0^1 u^{-\frac{pd}{2}}\left(\int_0^1 
	e^{-\frac{pv^2}{2\kappa u}}\,\dd v\right)^d\,\dd u = \frac{1}{(2\pi\kappa)^{\frac{pd}{2}} 
		} \int_0^1 u^{\frac{d(1-p)}{2}} \left( \int_0^{\frac{1}{\sqrt{u}}} 
	e^{-\frac{py^2}{2\kappa}} \, \dd y\right)^d \,\dd u\\
	&\geq\frac{1}{(2\kappa \pi)^{\frac{pd}{2}}} \left( \int_0^{1} 
	e^{-\frac{py^2}{2\kappa}} \, \dd y\right)^d \int_0^1 u^{\frac{d(1-p)}{2}} \,\dd u,
\end{align*}
which is finite if and only if $p<1+2/d$. So we conclude that 
\[ \E[|Y(t,x)|^{1+\frac{2}{d}}]=\infty \]
for all $(t,x)\in(1,\infty)\times\bbr^d$, and, in fact for all $t > 0$. It is not 
surprising that this holds in a much more general setting. The following results also 
answers an open problem posed in \cite[Remark~1.5]{Balan16}. Its proof will be given after 
the proof of Theorem~\ref{intfronts}.

\bthm[Non-existence of high moments]\label{nomoments}
Consider the situation described in Proposition~\ref{thm:moments} and assume that $\la\not\equiv 0$. Furthermore, suppose that there exists $(t_0,x_0)\in(0,\infty)\times\R^d$ such that 
\beq\label{sigma-not-zero} \si(Y_0(t_0,x_0))\neq 0, \eeq
where $Y_0$ is defined in \eqref{Y0}. If $Y$ denotes the unique mild solution to \eqref{SHE}, then 
\beq\label{infmom}\sup_{(t,x)\in[0,T]\times\bbr^d} \bbe[|Y(t,x)|^{1+\frac{2}{d}}]= +\infty \eeq
for all $T>t_0$.
\ethm

\brem The arguments presented in \cite{Bertini95} linking the notion of weak intermittency as defined in Definition~\ref{def:interm} with physical intermittency remain valid even if $\un\ga(p)=\infty$ for large values of $p$, provided we have $\un\ga(p)\uparrow\infty$ for $p\uparrow p_\mathrm{max}=\inf\{p>0\colon \un\ga(p)=\infty\}\leq 1+2/d$. Under mild assumptions, this is indeed the case as we will see in Theorem~\ref{asymptotics-1}.
\erem

\subsection{The martingale case}

In this subsection, we assume that $\La$ has mean zero, that is, $b=0$. As in the 
Gaussian case, we cannot hope for weak intermittency of order $1$ in general. This is a 
consequence of the following comparison principle for the stochastic heat equation driven 
by a nonnegative pure-jump Lévy noise, whose proof we postpone to the end of 
Section~\ref{Proof3}. The Gaussian analogue was established in  
\cite[Theorem~3.1]{Mueller91}.

\bthm[Comparison principle]\label{compprin} Let $\si$ be a non-decreasing Lipschitz function and 
$\La$ be a Lévy noise as in \eqref{eq:Lambda} with $b\in\R$, $\rho=0$ and $\la$ 
satisfying $\la((-\infty,0])=0$ and $m_\la(p)<\infty$ for some $p\in[1,1+2/d)$. Assume 
that $f_1\geq f_2 \geq0$ are two bounded measurable initial conditions, and $Y_1$ and 
$Y_2$ the corresponding mild solutions to \eqref{SHE}. There exist modifications of $Y_1$ 
and $Y_2$ such that, with probability $1$, we have $Y_1(t,x)\geq Y_2(t,x)$ for all 
$(t,x)\in[0,\infty)\times\R^d$.

In particular, if we have in addition that $f$ is a bounded nonnegative function and $0\leq \si(x)\leq Lx$ for some $L>0$, then the mild solution $Y$ to 
\eqref{SHE} has a nonnegative modification with
\beq\label{mean-calc} e^{(b\wedge 0)L t} \int_{\R^d} g(t,x-y)f(y)\,\dd y\leq  \bbe[|Y(t,x)|] = \bbe[Y(t,x)] \leq e^{(b\vee 0)L t} \int_{\R^d} g(t,x-y)f(y)\,\dd y \eeq
for all $(t,x)\in(0,\infty)\times\R^d$.
So if $b=0$, we have
$\overline \gamma(1) = 0$ if $f$ is strictly positive on a set of positive Lebesgue measure;
$\underline \gamma(1) = 0$ if $\inf_{x \in \R^d} f(x) > 0$;
$\overline \lambda(1) =0$ if $f(x) = O(e^{-c |x|})$ for some $c >0$;
and $\underline \lambda(1) = 0$ by definition.
\ethm

Thus, we are left to consider exponents in the region $p\in(1,1+2/d)$. In dimension 
$1$, we can use Itô's isometry to calculate second moments, and there are 
essentially no differences to the estimates (or exact formulae) obtained in the Gaussian 
case (\cite{Chen15, Conus12, Foondun09}). However, for $d\geq2$, we cannot use Itô's 
isometry because $p$ is strictly between $1$ and $2$. Instead, our main tool for proving 
intermittency  in the regime $p<2$ are the following moment lower bounds for stochastic 
integrals with respect to compensated Poisson random measures, which are of independent 
interest and complement existing sharp (but for our purposes not feasible) estimates in 
the literature (see \cite{MR}).

\blem\label{Poi-ineq} 
Let $(\calf_t)_{t\geq0}$ be a filtration on the underlying probability space and $N$ be an $(\calf_t)_{t\geq0}$-Poisson random measure on 
$[0,\infty)\times E$, where $E$ is a Polish space. Further suppose that $m$ denotes the intensity measure of $N$, 
and  $H\colon \Om\times[0,\infty)\times E\to\bbr$ 
is an $(\calf_t)_{t\geq0}$-predictable process such that the process
\[ t\mapsto \int_0^t \int_E H(s,x)\,\tilde N(\dd s,\dd x) \]
is a well-defined $(\calf_t)_{t\geq0}$-local martingale, where $\tilde N(\dd t,\dd x) := N(\dd t,\dd x)-m(\dd t,\dd x)$ is the 
compensation of $N$.

Then there exists for every $p\in(1,2]$ a 
constant $C_p>0$ that is independent of $H$ and $m$ such that 
\beq\label{poiest} 
\bbe\left[ \left|\iint_{[0,\infty)\times E} H(t,x)\,\tilde N(\dd t,\dd 
x)\right|^p \right] \geq C_p \frac{\iint_{[0,\infty)\times E} 
\bbe[|H(t,x)|^p] \,m(\dd t,\dd x)}{(1\vee m([0,\infty)\times E))^{1-\frac{p}{2}}},
\eeq
where $\infty/\infty:=0$.
In particular, if the right-hand side of \eqref{poiest} is infinite, then also the left-hand side of \eqref{poiest} is infinite. Furthermore, for every $p^\prime\in(1,2]$, the constants $C_p$ can be chosen to be bounded away from $0$ for $p\in[p^\prime,2]$.
\elem 

We are now ready to state the intermittency lower bounds for \eqref{SHE} that complement 
the corresponding upper bounds in Theorem~\ref{upper-bounds}. We start with non-vanishing 
initial data.

\bthm[Intermittency lower bounds -- I] \label{interm}
Let $Y$ be the solution to \eqref{eq:heat} constructed under the assumptions of Proposition~\ref{thm:moments}. 
Additionally assume that 
\begin{equation} \label{eq:sigma-inter}
L_f := \inf_{x\in\bbr^d} f(x) >0\quad\text{and}\quad 
L_\si:= \inf_{x \in \R \setminus\{0\} } \frac{|\sigma(x)|}{|x|} >0,
\end{equation}
and that $\La$ has the properties
\beq\label{condLa} b=0,\quad \la\not\equiv0\quad \text{and}\quad \int_\R |z|^{1+\frac{2}{d}}\bone_{\{|z|>1\}}\,\la(\dd z)<\infty. \eeq
Then the following statements are valid.
\benu
\item There exists a value $p_0 = p_0(\La,\kappa,\si)\in[1,1+2/d)$ such that we have
$\un\ga(p) >0$ 
for all exponents $p_0< p<1+2/d$.
\item For given $p\in(1,1+2/d)$, there exists $\kappa_0=\kappa_0(\La,p,\si)\in(0,\infty]$ such that  $\un\ga(p)>0$ for all diffusion constants $0<\kappa<\kappa_0$.
\item Given $p\in(1,1+2/d)$ and $\kappa>0$, there exists $L_0=L_0(\La,p,\kappa)\in[0,\infty)$ such that $\un\ga(p)>0$ if $\si$ has the property $L_\si>L_0$.
\item In dimension $d=1$, we can take $p_0=1$, $\kappa_0=\infty$ and $L_0=0$. 
\eenu
\ethm

To paraphrase, under the assumptions of Theorem~\ref{interm}, we have weak intermittency 
of order $p$ for every $p\in(1,3)$ in dimension $1$, while for higher dimensions we have 
this if $p$ is close enough to $1+2/d$, or $\kappa$ is small enough, or the size of $\si$ 
(or equivalently, the noise intensity) is large enough. It remains an open question whether in 
dimension $d\geq2$, we always have intermittency of all orders $p\in(1,1+2/d)$. Also, in 
contrast to the jump case where we have an affirmative answer, it seems to be open whether 
the solution to \eqref{SHE} in $d=1$ with Gaussian noise is weakly intermittent of order 
$p\in(1,2)$. 

For decaying initial condition, we have the following counterpart for the indices $\un\la(p)$.

\bthm[Intermittency lower bounds -- II]\label{intfronts} 
Let $Y$ be the solution to \eqref{eq:heat} constructed in Proposition~\ref{thm:moments}. 
Further assume that $c>0$, $L_\si>0$ (as defined in \eqref{eq:sigma-inter}), $\si(0)=0$, 
 that $f$ is nonnegative 
and strictly positive on a set of positive Lebesgue measure, $f(x) = O(e^{-c |x|})$ 
as $|x| \to \infty$, and that $\La$ satisfies \eqref{condLa}. 
\benu
\item There exists a value $p_1=p_1(\La,\kappa,\si)\in[1,1+2/d)$ such that 
$\un\la(p)>0$
for all $p\in(p_1,1+2/d)$. 
\item Given $p\in(1,1+2/d)$, there exists $\kappa_1=\kappa_1(\La,p,\si)\in(0,\infty]$ such that $\un\la(p)>0$ for all $0<\kappa<\kappa_1$.
\item Given $p\in(1,1+2/d)$ and $\kappa>0$, there exists $L_1=L_1(\La,p,\kappa)\in[0,\infty)$ such that $\un\la(p)>0$ for all $\si$ satisfying $L_\si>L_1$.
\item In $d=1$, we can take $p_1=1$, $\kappa_1=\infty$ and $L_1=0$.
\eenu
\ethm

\brem If $d=1$, $m_\la(2)<\infty$ and we consider the indices $\ov\ga(2)$, $\un\ga(2)$, $\un\la(2)$ and $\ov\la(2)$, there is -- thanks to Itô's isometry -- absolutely no difference between a Lévy and a Gaussian noise if we replace $\si$ by $\sqrt{v}\si$ where $v=\rho^2 + m_\la(2)^2$ is the variance of $\La$. For example, the explicit formulae derived in \cite{Chen15} immediately extend to the Lévy case.
\erem

\brem In \cite{Chen15}, the authors consider the stochastic heat equation with a measure-valued (e.g., a Dirac delta) initial condition. Their proof for the existence and uniqueness of solutions can be adapted to the Lévy setting by replacing $L^2$-estimates with $L^p$-type estimates from the BDG inequalities. Furthermore, since the heat operator smooths out a rough initial condition immediately, the intermittency properties of the solution will only depend on its decay and support properties. For example, Theorem~\ref{intfronts} as well as the Theorems~\ref{pos-mean}(2), \ref{interm-neg} and \ref{asymptotics-1}(2) continue to hold for the solution with a Dirac delta initial condition. 
\erem

\brem\label{origin} The intermittency of \eqref{SHE} with Gaussian noise is analytically 
due to the non-integrable tails of $g^2$ at $t=+\infty$ (see \cite{Conus12, Foondun09}). 
Translated into the picture of physical intermittency, this suggests that peaks in the 
past remain ``visible'' for a long time, and finally add up to new peaks. In the Lévy 
case, our proofs hint at the same phenomenon in dimension $1$ for the intermittency 
islands of low order (i.e., $p$ close to $1$). However, regardless of dimension, peaks of 
orders close to $1+2/d$, which are the dominating ones from a macroscopic level,  arise 
from  the singularity of the heat kernel at small times (this is further confirmed in the 
asymptotics we derive in Theorem~\ref{asymptotics-1}). It seems that high-order 
intermittency islands immediately trigger the formation of similar (or even larger) 
islands, leading to ``clusterings'' of peaks. It would be interesting for future research 
to specify and prove these heuristics.
\erem

\subsection{Noise with positive or negative drift}

In this section we consider the intermittency problem for \eqref{SHE} when the noise $\La$ has a non-zero mean. If $\La$ has a positive mean, that is, if $b>0$, then under natural assumptions, the solution to \eqref{SHE} is even weakly intermittent of order $1$ (and hence also of all orders $p\in[1,1+2/d)$).

\bthm[Intermittency for noises with positive drift]\label{pos-mean}
Suppose that $Y$ is the solution to \eqref{SHE} constructed in Proposition~\ref{thm:moments} and assume that $\si$ is a nonnegative Lipschitz continuous function with $L_\si>0$ (as defined in \eqref{eq:sigma-inter}). Furthermore, if $c=0$, suppose that $L_f$, as defined in \eqref{eq:sigma-inter}, is strictly positive, while for $c>0$, suppose that $f$ is nonnegative and strictly positive on a set of positive Lebesgue measure. If $b>0$, the following statements are valid.
\benu
\item If $c=0$, then $\un\ga(1)>0$.
\item If $c>0$, then $\un\la(1)>0$.
\eenu
\ethm

If $\La$ has a negative drift, we restrict ourselves to the \emph{parabolic Anderson model} where $\si$ is given by \eqref{PAM}. In this case, we can reformulate \eqref{SHE} as an equation driven by the martingale part of $\La$ only. In fact, decomposing $\La(\dd t,\dd x) = b\,\dd t\,\dd x + M(\dd t,\dd x)$, equation \eqref{SHE} can be written in the form
\beq\label{SHE-b}
\begin{split}
	\partial_t Y(t,x) &= \frac{\kappa}{2}\Delta Y(t,x) +b\si_0 Y(t,x)+ \si_0Y(t,x)\dot M (t,x),\quad 
	(t,x)\in(0,\infty)\times\bbr^d,\\
	Y(0,\cdot)&= f.
\end{split}
\eeq 
This is the $d$-dimensional \emph{stochastic cable equation} driven by the zero-mean 
Lévy space--time white noise $\dot M$. In a similar form, it has been studied in 
\cite{Walsh86} for Gaussian driving noise in dimension $d=1$. Its mild form is the same as 
in \eqref{eq:heat} but with $g$ replaced by
\[ g^\prime(t,x) = g(t,x)e^{b\si_0 t},\quad (t,x)\in(0,\infty)\times\R^d. \]

\bprop\label{cable} Under the assumptions of Proposition~\ref{thm:moments}, there exists 
$\beta_1>0$ such that \eqref{SHE-b} has a unique mild solution $Y$ satisfying 
$\|Y\|_{p,\beta,c}<\infty$ for all $\beta\geq\beta_1$. Furthermore, it is a modification 
of the unique mild solution to \eqref{SHE} constructed in Proposition~\ref{thm:moments}.
\eprop

We omit the proof since the existence and uniqueness result follows exactly as in the 
proof for Proposition~\ref{thm:moments}. Moreover, the second statement holds because weak 
and mild solutions are equivalent in our present setting: The proof is the same as in
\cite[Theorem~3.2]{Walsh86} for Gaussian $M$ and $d=1$. 

\bthm[Intermittency for noises with negative drift]\label{interm-neg}
Let $Y$ be the mild solution to \eqref{SHE} as in Proposition~\ref{thm:moments}. Suppose that $b<0$, $m_\la(1+2/d)<\infty$ and that $\si$ is given by \eqref{PAM} with $\si_0>0$. If $c=0$, also assume that $L_f>0$, and if $c>0$, that $f$ is nonnegative and positive on a set of positive Lebesgue measure.
\benu
	\item If $\la\not\equiv0$, Theorem~\ref{interm}(1)--(3) and Theorem~\ref{intfronts}(1)--(3) continue to hold.
	\item Let a value $p\in (1,1+2/d)$ be given, with the restriction $p\geq2$ if $\rho\neq0$. Whenever $\kappa$ or $|b|$ is large enough, or $\si_0$ is small enough (each time keeping the other two variables fixed), we have $\un\ga(p)\leq \ov\ga(p)<0$ and $\un\la(p)=\ov\la(p)=0$.
\eenu
\ethm

\section{Asymptotics of intermittency exponents}\label{Sect4}

As seen in the previous sections, the intermittency of the mild solution to \eqref{SHE} is stronger for higher values of $p$ or smaller values of $\kappa$. 
In this section, we investigate the limiting behavior of $\ov\ga(p)$, $\un\ga(p)$, $\ov\la(p)$ and $\un\la(p)$ as
\[ p\to 1+\frac{2}{d}\quad\text{and}\quad \kappa \to 0.
\]
In \eqref{kappato0} and \eqref{kappato0-2} below, one should keep in mind that, although not explicitly indicated in the notation, the indices $\ov\ga(p)$ etc.\ also depend on $\kappa$.
\bthm[Asymptotics of intermittency exponents]\label{asymptotics-1} Consider a noise $\La$ with non-zero Lévy measure $\la$.
\benu \item Let $c=0$ and grant the assumptions of Theorem~\ref{interm}, Theorem~\ref{pos-mean} or Theorem~\ref{interm-neg} depending on whether $\La$ has mean $b=0$, $b>0$ or $b<0$. If $b>0$ or $b<0$, we also impose that $\si$ is of the form \eqref{PAM}. Then we have
\beq\label{p-upperbound}  
\lim_{p\to 1+\frac{2}{d}} 
\frac{1+\frac{2}{d}-p }{\left|\log\left( 1+\frac{2}{d}-p  \right)\right|}
\log \un\ga(p) =  \lim_{p\to 1+\frac{2}{d}} 
\frac{1+\frac{2}{d}-p }{\left|\log\left(1+\frac{2}{d}-p  \right)\right|}
\log \ov\ga(p)=\frac{2}{d},
\eeq
\beq\label{kappato0} 0<\liminf_{\kappa\to0} \kappa^{\frac{p-1}{1+2/d-p}} \un\ga(p) \leq \limsup_{\kappa\to 0} \kappa^{\frac{p-1}{1+2/d-p}} \ov\ga(p)<\infty. \eeq
\item Let $c>0$ and grant the assumptions of Theorem~\ref{intfronts}, Theorem~\ref{pos-mean} or Theorem~\ref{interm-neg} depending on whether $\La$ has mean $b=0$, $b>0$ or $b<0$. If $b>0$ or $b<0$, we also impose that $\si$ is of the form \eqref{PAM}. Then we have
\beq\label{p-upperbound-2} \frac{1}{d}\leq \liminf_{p\to 1+\frac{2}{d}} 
\frac{1+\frac{2}{d}-p }{\left|\log\left( 1+\frac{2}{d}-p  \right) \right|}
\log \un\la(p) \leq \limsup_{p\to 1+\frac{2}{d}} 
\frac{ 1+\frac{2}{d}-p }{\left|\log\left( 1+\frac{2}{d}-p  \right)\right|} \log 
\ov\la(p) \leq \frac{2}{d}.\eeq
If in addition the initial condition decays superexponentially in the sense that 
$|f(x)|=O(e^{-c|x|})$ as $|x|\to\infty$ for every $c\geq0$, then 
\beq\label{kappato0-2} 0<\liminf_{\kappa\to0} \kappa^{-\frac{1+1/d-p}{1+2/d-p}} \un\la(p) \leq \limsup_{\kappa\to 0} \kappa^{-\frac{1+1/d-p}{1+2/d-p}} \ov\la(p)<\infty. \eeq
\eenu
\ethm

\brem
\benu
\item 
Equation~\eqref{p-upperbound} asserts that the moment Lyapunov exponents $\un\ga(p)$ and $\ov\ga(p)$, which determine the exponential rates at which $\bbe[|Y(t,x)|^p]$ grows for $t\to\infty$, themselves increase at a superexponential speed as $p$ approaches $1+2/d$. This is much faster than in the Gaussian case, where for the PAM \eqref{PAM} in $d=1$ with constant $f$ \cite[Theorem~2.6]{Bertini95} and \cite[Theorem~6.4]{Khos} showed that the Lyapunov exponents have a cubic growth as $n\to\infty$:
\beq\label{PAM-Gauss} \un\ga(n)=\ov\ga(n)=\frac{\si_0^4}{4! \kappa}n(n^2-1),\quad n\in\N. \eeq
We conclude that the intermittent behavior of the stochastic heat equation with jumps is much stronger than with Gaussian noise. 
\item  Similarly, \eqref{p-upperbound-2} states that the velocity at which $p$th order intermittency peaks propagate in space grows superexponentially when $p\to1+2/d$. Again, this is on a much faster scale than in the Gaussian case, where the indices $\un\la(p)$ and $\ov\la(p)$ typically only increase linearly in $p$: see \cite[Proposition~3.11]{Hu16} where for the PAM \eqref{PAM} in $d=1$ with compactly supported initial data $f$, the authors showed that
\beq\label{PAM-Gauss-lambda} 0<\liminf_{n\to\infty} \frac{\un\la(n)}{n} \leq \limsup_{n\to\infty} \frac{\ov\la(n)}{n}<\infty.  \eeq
We also remark that in the jump case, the asymptotics of the exponents $\un\ga(p)$ and $\ov\ga(p)$ as $p\to1+2/d$ are similar to the exponents $\un\la(p)$ and $\ov\la(p)$, in contrast to the Gaussian case, cf.\ \eqref{PAM-Gauss} and \eqref{PAM-Gauss-lambda}.
\item Regarding the asymptotics for $\kappa>0$, a notable difference between jump and Gaussian noise is that in the former case, the rate at which $\un\ga(p)$ and $\ov\ga(p)$ increases as $\kappa\to0$ explicitly depends on $p$, whereas in the latter case, at least for $p\in\N$, it typically does not, see \eqref{PAM-Gauss}.
\item Another interesting observation is that for jump noises, the asymptotics of 
$\un\la(p)$ and $\ov\la(p)$ for $\kappa\to0$ exhibit a phase transition at $p=1+1/d$. If 
$p\in(1,1+1/d)$, they decrease like $\kappa^{(1+1/d-p)/(1+2/d-p)}$, if $p=1+1/d$, they are 
bounded away from zero and infinity in $\kappa$, and for $p\in(1+1/d,1+2/d)$, they 
increase like $\kappa^{-(p-1+1/d)/(1+2/d-p)}$. Intuitively speaking, this is because for 
small $\kappa$ there are two effects that counteract each other: On the one hand, a small 
diffusion constant reduces the speed at which the initial mass at the origin can spread. 
On the other hand, if $\kappa$ is small, once an intermittency peak is built up, it takes 
longer for the Laplace operator to smooth it out, which facilitates the development and 
transmission of further peaks. Thus, for small values of $p$, the first effect is 
dominant, while for large values of $p$, it is the second effect that wins. In the 
Gaussian case, the behavior is again different. Here for any $p\in[2,\infty)$, we have
\beq 0<\liminf_{\kappa\to0} \un\la(p)\leq \limsup_{\kappa\to0} \ov\la(p)<\infty. \eeq
The lower bound follows from \cite[Theorem~1.3]{Conus12} together with the fact that $\un\la(2) \leq \un\la(p)$ for all $p\geq2$, while the upper bound follows as in the proof of Theorem~\ref{asymptotics-1} from the formula \eqref{Cbetagamma}.
\eenu
\erem

\section{Proofs}
\subsection{Proofs for Section~\ref{Sect2}}
\blem\label{kernel-int} 
Define $g_{\beta,c}(t,x) := g(t,x)e^{-\beta t+c|x|}$ for $(t,x)\in(0,\infty)\times\bbr^d$. If $0<p<1+2/d$, $c\geq0$ and $\beta>\kappa c^2 d/2$, then
\[ 
\int_0^\infty \int_{\bbr^d} g^p_{\beta,c}(t,x) \,\dd t\,\dd x \leq 
\frac{2^{\frac{d}{2}(3-p)}\Gamma(1-\frac{d}{2}(p-1))}{p^{1+d(1-\frac{p}{2})}(\pi\kappa)^{\frac{d}{2}(p-1)} (\beta-\frac12 \kappa c^2 
	d)^{1-\frac{d}{2}(p-1)}},
\]
where $\Ga$ denotes the gamma function $\Ga(x)=\int_0^\infty t^{x-1}e^{-t}\,\dd t$.
\elem

\bpr If $\beta>\kappa c^2 d/2$, then
\begin{align*}
	\int_0^\infty\int_{\R^d} g^p_{\beta,c}(t,x)\,\dd t\,\dd x &=\int_0^\infty 
	\frac{e^{-p\beta t}}{p^{\frac{d}{2}}(2\pi \kappa t)^{\frac{d}{2}(p-1)}} \int_{\bbr^d} 
	\frac{e^{-\frac{p}{2\kappa t}|x|^2}}{(2\pi \kappa t/p)^{\frac{d}{2}}} e^{p c |x|}\,\dd x\,\dd t\\
	&\leq \int_0^\infty \frac{e^{-p\beta t}}{p^{\frac{d}{2}}(2\pi\kappa t)^{\frac{d}{2}(p-1)}} 
	\left(\int_{\bbr} \frac{e^{-\frac{p}{2\kappa t}|x|^2}}{(2\pi\kappa t/p)^{\frac12}} e^{pc |x|}\,\dd 
	x\right)^d\,\dd t\\
	&\leq \int_0^\infty \frac{2^d e^{-p\beta t}}{p^{\frac{d}{2}}(2\pi\kappa t)^{\frac{d}{2}(p-1)}} 
	\left(\int_{\bbr} \frac{e^{-\frac{p}{2\kappa t}|x|^2}}{(2\pi\kappa t/p)^{\frac12}} e^{pc x}\,\dd 
	x\right)^d\,\dd t\\
	&= \int_0^\infty \frac{2^d e^{-p\beta t}}{p^{\frac{d}{2}}(2\pi\kappa t)^{\frac{d}{2}(p-1)}} e^{\frac{1}{2}d \kappa
		p c^2 t}\,\dd t \\
	&= 
	\frac{2^{\frac{d}{2}(3-p)}\Gamma(1-\frac{d}{2}(p-1))}{p^{1+d(1-\frac{p}{2})}(\pi\kappa)^{\frac{d}{2}(p-1)} (\beta-\frac12 \kappa c^2 
		d)^{1-\frac{d}{2}(p-1)}}.
\end{align*}
\epr 

\begin{proof}[\bff{\emph{Proof of Proposition~\ref{prop:stochYoung-d}}}]
	We use the triangle inequality to split
	\[
	\begin{split}
	\left\| (g \circledast \Phi)(t,x) \right\|_p & \leq  
	|\rho| \left\| \int_0^t \int_{\R^d} g(t-s, x-y) \Phi(s,y)\, W( \dd s, \dd y) \right\|_p \\
	& \phantom{\leq} + 
	\left\| \int_0^t \int_{\R^d} \int_{\R} g(t-s, x-y) 
	\Phi(s,y) z\, (\mu - \nu) ( \dd s, \dd y, \dd z) \right\|_p \\
	& \phantom{\leq} + 
	|b| \left\| \int_0^t \int_{\R^d} g(t-s, x-y) \Phi(s,y)\,  \dd s \,\dd y \right\|_p \\
	& =: I_1(t,x) + I_2(t,x) + I_3(t,x)
	\end{split}
	\]
	into a Gaussian, a pure-jump and a drift part. Recall that $I_1$ vanishes for $d\geq2$. For $d=1$ and $p\in[2,3)$, 
	we have from the BDG inequalities (see \cite[Theorem~VII.92]{Dellacherie82}) together with Minkowski's 
	integral inequality that 
	\begin{align*} 
		e^{-\beta t + c|x|}I_1(t,x) 
		&\leq 
		|\rho|C_pe^{-\beta t + c|x|}\left( \displaystyle\int_0^t \int_\bbr 
		g^2(t-s,x-y)\|\Phi(s,y)\|_p^2 \,\dd s\,\dd y \right)^{\frac12} \\
		&\leq |\rho| C_p \|\Phi\|_{p,\beta,c} \left( \int_0^t\int_\R 
		g^2(t-s,x-y) e^{-2\beta(t-s)+2 c (|x|-|y|)}\,\dd s\,\dd y  \right)^{\frac12}\\
		&\leq |\rho|C_p \|\Phi\|_{p,\beta,c} \left( \int_0^\infty\int_\R 
		g^2_{\beta,c}(s,y)\,\dd s\,\dd y  \right)^{\frac12}.
	\end{align*}
	So we deduce from Lemma~\ref{kernel-int} 
	that 
	\beq
	\label{eq:I1} 
	\sup_{(t,x)\in(0,\infty)\times\bbr} e^{-\beta t + c|x|}I_1(t,x) \leq 
	C_p|\rho| \frac{1}{(2\kappa(\beta- \frac12 \kappa c^2))^{\frac14}}\|\Phi\|_{p,\beta,c}.
	\eeq

	In order to estimate $I_3$ we only need Minkowski's integral inequality and Lemma~\ref{kernel-int} to obtain
	\begin{equation} \label{eq:I3}
		\begin{split}
			I_3(t,x) &\leq |b| \int_0^t \int_{\R^d} g(t-s, x-y) 
			\| \Phi(s,y)\|_p\, \dd s\, \dd y\\
			& \leq |b| e^{\beta t- c|x|}\|\Phi\|_{p,\beta,c}  \int_0^t\int_{\R^d} g(t-s,x-y) e^{-\beta (t-s) +c(|x|-|y|)} \,\dd s\,\dd y\\
			& \leq |b| e^{\beta t-c|x|}\|\Phi\|_{p,\beta,c}  \int_0^t\int_{\R^d} g(t-s,x-y) e^{-\beta (t-s) +c(|x-y|)} \,\dd s\,\dd y\\  
			&\leq \frac{2^d |b|}{\beta-\frac12 \kappa c^2 d} e^{\beta t-c|x|} 
			\|\Phi\|_{p,\beta,c}.
		\end{split}
	\end{equation}
	
	We turn to the estimation of $I_2$. If $p\leq 2$, we use the BDG
	inequality  to deduce
	\begin{align*}
		I_2(t,x)^p &\leq C^p_p \, \bbe\left[\left( \int_0^t\int_{\bbr^d}\int_\bbr 
		|g(t-s,x-y)\Phi(s,y)z|^2\,\mu(\dd s,\dd y,\dd z) \right)^{\frac{p}{2}}\right]\\
		&\leq C^p_p  \int_0^t\int_{\bbr^d}\int_\bbr g^p(t-s,x-y) \|\Phi(s,y)\|_p^p 
		|z|^p\,\nu(\dd s,\dd y,\dd z)\\
		&\leq  C^p_p (m_\la(p))^p  e^{p\beta t - p c|x|} 
		\frac{2^{\frac{d}{2}(3-p)}\Gamma(1-\frac{d}{2}(p-1))}{p^{1+d(1-\frac{p}{2})}(\pi\kappa)^{\frac{d}{2}(p-1)} (\beta-\frac12 \kappa c^2 
			d)^{1-\frac{d}{2}(p-1)}} \|\Phi\|^p_{p,\beta,c}.
	\end{align*}
	At the second inequality we used that $(\sum_{i=1}^\infty a_i)^r \leq \sum_{i=1}^\infty 
	a_i^r$ for any $r \in [0,1]$ and nonnegative numbers $(a_i)_{i\in\N}$.
	If $d=1$ and $2<p<3$, we use 
	\cite[Theorem~1]{MR} with $\alpha = 2$ to obtain
	\begin{equation} \label{eq:I2-1}
		\begin{split}
			I_2(t,x)^p & \leq C^p_p \Bigg(  \E\left[ \left( \int_0^t \int_\R \int_\R |g(t-s, x-y) 
			\Phi(s,y)  z|^2 \,
			\nu ( \dd s, \dd y, \dd z) \right)^{\frac{p}{2}}\right] \\
			&\quad \ +  \int_0^t \int_\R \int_\R g^p(t-s, x-y) \|\Phi(s,y)\|_p^p |z|^p \, 
			\nu ( \dd s, \dd y, \dd z) \Bigg).
		\end{split}
	\end{equation}
	For the first term, again by Minkowski's integral inequality and Lemma~\ref{kernel-int}, we have
	\begin{equation} \label{eq:I2-2}
		\left( \E \left[ \left(\int_0^t \int_\R \int_\R |g(t-s, x-y) \Phi(s,y)  z|^2 \,
		\nu ( \dd s, \dd y, \dd z) \right)^{\frac{p}{2}}\right]\right)^{\frac{1}{p}} \leq \frac{m_\lambda(2) e^{\beta t - c|x|}}{(2\kappa(\beta- \frac12 \kappa c^2))^{\frac14}}
		\|\Phi\|_{p,\beta,c},
	\end{equation}
	while for the second term, 
	\begin{equation} \label{eq:I2-3}
		\left(\int_0^t \int_\R \int_\R g^p(t-s, x-y) \|\Phi(s,y)\|_p^p |z|^p \, 
		\nu ( \dd s, \dd y, \dd z) \right)^{\frac{1}{p}}  
		\leq  	\frac{2^{\frac{3-p}{2p}}\Gamma(\frac{3-p}{2})^{\frac{1}{p}}m_\lambda(p) e^{\beta t - c|x|} }{p^{\frac{4-p}{2p}}(\pi\kappa)^{\frac{p-1}{2p}} (\beta-\frac12 \kappa c^2 
			)^\frac{3-p}{2p}} \|\Phi\|_{p,\beta,c}.
	\end{equation}
	Substituting  \eqref{eq:I2-2} and \eqref{eq:I2-3} back into \eqref{eq:I2-1}, we obtain
	\begin{equation} \label{eq:I2}
		e^{-\beta t + c|x|}I_2(t,x) \leq C_p \left( 
		\frac{m_\lambda(2) }{(2\kappa(\beta-\frac12 \kappa c^2))^{\frac14}} + 
		\frac{2^{\frac{3-p}{2p}}\Gamma(\frac{3-p}{2})^{\frac{1}{p}}m_\lambda(p)  }{p^{\frac{4-p}{2p}}(\pi\kappa)^{\frac{p-1}{2p}} (\beta-\frac12 \kappa c^2 
			)^\frac{3-p}{2p}} \right)
		\|\Phi\|_{p,\beta,c}.
	\end{equation}
	The statement now follows from inequalities \eqref{eq:I1}, \eqref{eq:I3} and 
\eqref{eq:I2}. Finally, since $C_p$ comes from BDG inequalities, it remains bounded on 
$[1+\epsilon, 1+2/p)$. 
\end{proof}

\begin{proof}[\emph{\bff{Proof of Proposition~\ref{thm:moments}}}]
	The proof combines Proposition~\ref{prop:stochYoung-d} with arguments in \cite[Theorem~1.1]{Conus12} (see also \cite[Theorem~8.1]{Khos}).
	
 As usual, we consider the Picard iteration sequence $Y^{(0)} = Y_0$ and 
	\[
	Y^{(n)}(t,x) =  Y_0(t,x) + \int_0^t \int_{\R^d} g(t-s, x-y) \sigma(Y^{(n-1)}(s,y)) \,
	\Lambda(\dd s, \dd y)
	\]
	for $n\in\N$, and define $u^{(n)}=Y^{(n)}-Y^{(n-1)}$. 	After possibly enlarging the value of $L$, we can assume that $|\si(x)|\leq L(1+|x|)$ for all $x\in\bbr$. Now
	let us choose $\beta_0>\frac12\kappa c^2 d$ large enough such that the factor $C_{\beta,c}(\kappa,p)$ in front of 
	$\|\Phi\|_{p,\beta,c}$ on the right-hand side of \eqref{eq:stochYoung-ineq-d} satisfies
	\beq\label{satisfies} C_{\beta,c}(\kappa, p)< \frac 1 L\quad \text{for all }\beta\geq\beta_0.\eeq
	Using the Lipschitz property of $\si$, we obtain for all $\beta\geq\beta_0$ and $n\in\bbn$ 
	as a consequence of Proposition~\ref{prop:stochYoung-d},
	\begin{align*}
		\|u^{(n)}\|_{p,\beta,c} &= \|g\circledast (\si(Y^{(n-1)})-\si(Y^{(n-2)})) 
		\|_{p,\beta,c} \leq C_{\beta,c}(\kappa,p) \| \si(Y^{(n-1)})-\si(Y^{(n-2)}) \|_{p, \beta, c}
		\\
		& \leq q \| u^{(n-1)}\|_{p,\beta,c} \leq \ldots\leq q^{n-1}\|u^{(1)}\|_{p,\beta,c}
	\end{align*}
	for some $q=q_c(\kappa,p)<1$. If $c=0$, the last term is less than or equal 
to 
	$Cq^{n}(1+\|Y_0\|_{p,\beta,c})$, while it is bounded by $Cq^{n}\|Y_0\|_{p,\beta,c}$ if 
	$c>0$ (and therefore $\si(0)=0$). Since $\beta\geq\beta_0>\frac12 \kappa c^2d$,
	\begin{align*} 
		\|Y_0\|_{p,\beta,c} &= \sup_{t \in (0,\infty)} \sup_{x \in \R^d} e^{-\beta t+c|x|} 
		\left|\int_{\bbr^d} g(t,x-y)f(y)\,\dd y\right|\\
		& \leq  \sup_{x\in\bbr^d} e^{c |x|}|f(x)| \sup_{t\in(0,\infty)} e^{-\beta t} 
		\int_{\bbr^d} g(t,x)e^{c|x|}\,\dd x \\
		& \leq C  \sup_{t\in(0,\infty)} e^{-\beta t} \left(\int_\bbr g(t,x) e^{c x}\,\dd 
		x\right)^d\\
		& = C\sup_{t\in(0,\infty)} e^{-(\beta -\frac12 \kappa c^2 d) t} <\infty,  
	\end{align*}
	it follows that $(Y^{(n)})_{n\in\bbn}$ is a Cauchy 
	sequence with respect to $\|\cdot\|_{p,\beta,c}$, converging in 
	$\|\cdot\|_{p,\beta,c}$ to some limit $Y$. That $Y$ satisfies \eqref{eq:heat} and is 
	unique up to modifications, follows as in \cite[Theorem~3.1]{Chong}.
\end{proof}

\bpr[\emph{\bff{Proof of Theorem~\ref{upper-bounds}}}] The first part follows immediately from $\|Y\|_{p,\beta,0}<\infty$ for 
$\beta\geq\beta_0$ with $\beta_0$ as in the proof of Proposition~\ref{thm:moments}.
Concerning the second part of the theorem, observe that $\|Y\|_{p,\beta,c}<\infty$ for 
$\beta \geq \beta_0$ implies $\bbe[|Y(t,x)|^p] \leq Ce^{\beta pt - c p|x|}$ for all 
$t>0$, $x \in \bbr^d$  and some finite constant $C>0$. Hence,
\[
\sup_{|x|\geq \al t}  \bbe[|Y(t,x)|^p] \leq Ce^{\beta pt - c p\al t},
\]
and therefore,
\beq\label{help}
\limsup_{t\to\infty} \frac{1}{t}\sup_{|x|\geq \al t} \log \bbe[|Y(t,x)|^p] < 0
\eeq
for all $\al>\beta_0/c$. 
\epr

\subsection{Proofs for Section~\ref{Sect3}}\label{Proof3}

\blem\label{poimom} 
If $X_\la$ has a Poisson distribution with parameter $\la$, then there exists for every 
$r>0$ a constant $C_r>0$ such that 
\[ \bbe[X^r_\la] \geq C_r \begin{cases}  \la^r&\text{for } \la>1,\\ \la&\text{for } \la\leq1.\end{cases} \]
\elem

\bpr 
Suppose that $(X_\la)_{\la\geq0}$ forms a standard Poisson process. The law of 
large numbers implies that $X_\la/\la \to 1$ a.s.\ as $\la\to\infty$. The convergence also 
takes place in $L^p$ for every $p\geq1$ because $\bbe[X_\la^n]$ is a polynomial in $\la$ 
of degree $n$ for every $n\in\bbn$ so that $\sup_{\la\geq1} \bbe[X_\la^n]/\la^n<\infty$. 
In particular, we obtain for every $r>0$ that $\bbe[X_\la^r]/\la^r \to 1$ as 
$\la\to\infty$, which implies the claim for $\la>1$. The bound for $\la\leq1$ follows 
from the definition of the expectation and 
$\bbp[X_\la^r=1]=\bbp[X_\la=1]=\la e^{-\la}\geq\la e^{-1}$.
\epr 

The following \emph{decoupling inequalities} can be found in 
\cite[Theorem~2.4.1]{Veraar06}. Because of its importance for proving 
Lemma~\ref{Poi-ineq}, and because the proof in the reference is given for processes with 
values in Banach spaces, we reproduce the proof in the real-valued setting for the 
reader's convenience.
In the following lemma, for notational ease, a random variable $\xi\colon \Omega \to 
\mathbb{R}$ is identified with its natural extension to the product space $\Omega \times 
\overline \Omega$, i.e., $\xi(\omega, \overline \omega) = \xi(\omega)$.

\blem\label{decoupling}
Consider two probability spaces $(\Om, \calf,\bbp)$ and $(\overline \Om,\overline 
\calf,\overline \bbp)$, each of them equipped with a discrete-time filtration 
$(\calf_i)_{i\geq0}$ and $(\overline \calf_i)_{i\geq0}$, respectively. Furthermore, let  
$(\xi_i)_{i\geq1}$ be a zero-mean $(\calf_i)_{i\geq1}$-adapted sequence such that 
$\xi_{i}$ is independent of $\calf_{i-1}$ under $\bbp$ for all $i\geq1$, and let 
$(\overline \xi_i)_{i\geq1}$ be a sequence with analogous properties on $(\overline  
\Om,\overline  \calf,\overline \bbp)$ and the same distribution as $(\xi_i)_{i\geq1}$. 
Finally, assume that $(H_i)_{i\geq1}$ is a sequence of random variables on 
$(\Om\times\overline \Om,\calf\otimes\overline \calf,\bbp\otimes\overline \bbp)$ such 
that 
$H_i$ is $\calf_{i-1}\otimes\overline \calf_{i-1}$-measurable for all $i\geq1$. Then for 
every $p\in(1,\infty)$ there exist constants $C_p,C'_p>0$ that are independent of 
$(\xi_i)_{i\geq1}$ and $(H_i)_{i\geq1}$ such that for every $N\in\bbn$,
\[ 
(C'_p)^{-1} \bbe\left[ \overline \bbe\left[ \left| \sum_{i=1}^N H_i \overline  \xi_i \right|^p \right]  \right]
\leq \bbe\left[ \overline \bbe\left[\left| \sum_{i=1}^N H_i \xi_i \right|^p \right]\right] \leq 
C_p\bbe\left[ \overline \bbe\left[\left| \sum_{i=1}^N H_i \overline  \xi_i \right|^p \right]\right]. 
\]
\elem

\bpr Define the random variables
\[ D_{2i-1}:=\frac12 (H_i \xi_i + H_i \overline  \xi_i),\quad D_{2i} := \frac12(H_i\xi_i - 
H_i 
\overline  \xi_i),\quad i=1,\ldots,N, \]
and a filtration $(\calg_i)_{i=0,\ldots,2N}$ by
\[ \calg_0:=\{\emptyset,\Om\},\quad \calg_{2i-1}:=\si(\calf_{i-1} \otimes \overline  \calf_{i-1},\xi_i+\overline  
\xi_i),\quad 
\calg_{2i} := \calf_i \otimes \overline \calf_i, \quad i=1,\ldots,N. \]
Obviously, $(D_i)_{i=1,\ldots,2N}$ is adapted to $(\calg_i)_{i=1,\ldots,2N}$. In 
addition, denoting by $\bbe\otimes \ov\bbe$ the expectation with respect to $\bbp\otimes 
\overline \bbp$, we have for all $i=1,\ldots,N$,
\begin{align*} \bbe\otimes \ov\bbe[ D_{2i+1} \mid \calg_{2i} ] &= \frac12 
H_{i+1}\bbe\otimes \ov\bbe[\xi_{i+1}+\overline  \xi_{i+1}] = 0, \\
\bbe\otimes \ov\bbe[D_{2i} \mid \calg_{2i-1} ] &= \frac12 H_i\bbe\otimes \ov\bbe[\xi_i-\overline  
\xi_i\mid 
\xi_{i}+\overline  \xi_i] = 0, \end{align*}
where the last identity holds because $\xi_i$ and $\overline  \xi_i$ are independent with 
the 
same distribution. It follows from \cite[Theorem~VII.1.1]{Shiryaev96} that the processes $(\sum_{i=1}^n D_i)_{n=0,\ldots,2N}$ and 
$(\sum_{i=1}^n (-1)^{i+1} D_i)_{n=0,\ldots,2N}$ are discrete-time local martingales with respect to 
$(\calg_i)_{i=0,\ldots,2N}$.

Observing that 
\[ \sum_{i=1}^N H_i\xi_i = \sum_{i=1}^{2N} D_i,\quad \sum_{i=1}^N H_i \overline  \xi_i  = 
\sum_{i=1}^{2N} (-1)^{i+1}D_i \]
by construction, the claim is a consequence of the 
classical BDG inequalities because the two discrete-time local martingales above can be canonically embedded into continuous-time local martingales with the same quadratic variation 
process.
\epr

\bpr[\bff{\emph{Proof of Lemma~\ref{Poi-ineq}}}] We first prove \eqref{poiest} for simple integrands of the form
\beq\label{simple} 
Z(\om,t,x)=\sum_{i,j=1}^{K} X_{ij}(\om)\bone_{(t_{i-1},t_{i}]\times 
	B_j}(t,x),\quad(\om,t,x)\in\Om\times[0,\infty)\times E,
\eeq
where $0\leq t_0 \leq \ldots \leq t_K<\infty$, $(B_j)_{j=1,\ldots,K}$ are 
 pairwise disjoint
Borel subsets of $E$, and $X_{ij}$ are $\calf_{t_{i-1}}$-measurable 
random variables for all $i,j=1,\ldots,K$. 

Using Lemma~\ref{decoupling}, we can assume without loss of generality that 
$Z$ is deterministic, that is, the variables $X_{ij}(\om)$ do not depend on $\omega$. 
To see this, 
define
\[
\begin{split}
& \xi_{ij}(\omega)  =   N((t_{i-1}, t_i] \times B_j)(\omega) 
- m((t_{i-1}, t_i] \times B_j), \\
& \ov\xi_{ij}(\ov\omega)  =   \ov N((t_{i-1}, t_i] \times B_j)(\ov\omega) 
- m((t_{i-1}, t_i] \times B_j), \\
& H_{ij}(\omega, \overline \omega ) =  X_{ij}(\omega),
\end{split}
\]
where $\ov N$ lives on a copy $(\ov\Om,\ov\calf,(\ov\calf_t)_{t\geq0},\ov\bbp)$ of the original probabilty space, with the same distribution as $N$. Since $\xi_{ij}$ is $\mathcal{F}_{t_i}$-measurable 
and $H_{ij}$ is $\calf_{t_{i-1}}\otimes \ov\calf_{t_{i-1}}$-measurable, 
Lemma~\ref{decoupling} applies and yields
\[
\begin{split}
& \E \left[\left| \sum_{i,j=1}^K X_{ij}(\omega) \left( N((t_{i-1}, t_i] \times B_j)(\om) - 
m((t_{i-1}, t_i] \times B_j) \right) \right|^p \right]  \\
& \geq (C_p')^{-1}   \E \left[ \overline\E\left[ \left| \sum_{i,j=1}^K X_{ij}(\omega) 
\left( \overline N((t_{i-1}, t_i] \times B_j)(\overline \omega) - m((t_{i-1}, t_i] 
\times B_j) \right) \right|^p\right] \right].  
\end{split}
\]
As $X_{ij}(\omega)$ does not depend on $\overline \omega$, 
it is indeed enough to prove \eqref{poiest}
for deterministic integrands.

By the BDG inequalities, there exists $C_p>0$ (which is bounded away from $0$ for 
$p>p^\prime$) such that
\begin{align*} 
	\bbe\left[ \left|\iint_{[0,\infty)\times E} Z(t,x)\,\tilde N(\dd t,\dd x)\right|^p \right] 
	&\geq C_p \bbe\left[ \left|\iint_{[0,\infty)\times E} Z^2(t,x)\, N(\dd t,\dd x)\right|^{\frac{p}{2}} 
	\right]\nonumber \\&= C_p\bbe\left[ \left|\sum_{i,j=1}^K X^2_{ij}N((t_{i-1},t_i] \times 
	B_j)\right|^{\frac{p}{2}}  \right].  
\end{align*}
Inequality \eqref{poiest} is shown for integrands of the form \eqref{simple} once we can 
show that 
\beq\label{poiest2} 
\bbe\left[ \left|\sum_{i=1}^K a_{i}N(A_i)\right|^{r} \right] \geq C \frac{\sum_{i=1}^K 
	a_i^r m(A_i)}{(1\vee m([0,\infty)\times E))^{1-r}} 
\eeq 
for all $a_i\in[0,\infty)$, pairwise disjoint $A_i\in\calb([0,\infty)\times E)$ and 
$r\in(1/2,1]$. By the tower property of conditional expectations, 
\begin{align*}
	\bbe\left[ \left( \sum_{i=1}^K a_iN(A_i) \right)^r \right] &= \sum_{n=1}^\infty \bbe\left[ 
	\left( \sum_{i=1}^K a_iN(A_i) \right)^r \,\Bigg\vert\, \sum_{i=1}^K N(A_i) = n \right] 
	\bbp\left[ \sum_{i=1}^K N(A_i) =n \right].
	\end{align*}
On the event $\sum_{i=1}^K N(A_i) = n$, at most $n$ summands in $\sum_{i=1}^K a_iN(A_i)$ are different from zero. Therefore, by rewriting $a_iN(A_i)$ as a sum $a_i+\ldots+a_i$ of $N(A_i)$ terms,  $\sum_{i=1}^K a_iN(A_i)$ becomes a sum of $\sum_{i=1}^K N(A_i)=n$ (possibly repeated) terms. Thus,
using the estimate $(\sum_{i=1}^n c_i)^r \geq n^{r-1}\sum_{i=1}^n c_i^r$ 
for nonnegative $c_1,\ldots,c_n$, we obtain
\begin{align*}
	\bbe\left[ \left( \sum_{i=1}^K a_iN(A_i) \right)^r \right]	&\geq \sum_{n=1}^\infty n^{r-1} \bbe\left[  \sum_{i=1}^K a_i^r N(A_i) \, \Bigg\vert\, 
	\sum_{i=1}^K N(A_i) = n \right]\bbp\left[ \sum_{i=1}^K N(A_i) =n \right]\\
	&=\sum_{i=1}^K  a_i^r \sum_{n=1}^\infty n^{r-1} \bbe\left[ N(A_i) \, \Bigg\vert\, 
	\sum_{i=1}^K N(A_i) = n  \right]\bbp\left[ \sum_{i=1}^K N(A_i) =n \right]\\
	&=\sum_{i=1}^K  a_i^r \sum_{n=1}^\infty n^{r} \frac{m(A_i)}{\sum_{j=1}^K m(A_j)}  
	\bbp\left[ \sum_{i=1}^K N(A_i) =n \right]\\
	&= \frac{\sum_{i=1}^K a_i^r m(A_i)}{\sum_{j=1}^K m(A_j)} \sum_{n=1}^\infty n^r \bbp\left[ 
	\sum_{i=1}^K N(A_i) =n \right]\\
	&=\frac{\sum_{i=1}^K a_i^r m(A_i)}{\sum_{j=1}^K m(A_j)}\bbe\left[ \left(\sum_{i=1}^K 
	N(A_i)\right)^r  \right].
\end{align*}
Since the constant $C_r$ in Lemma~\ref{poimom} can be taken independently of $r$ when 
$r \in (1/2,1]$, we derive
\[ 
\bbe\left[ \left( \sum_{i=1}^K a_i N(A_i) \right)^r \right]  \geq C \frac{\sum_{i=1}^K 
	a_i^r m(A_i)}{\left(1\vee\sum_{j=1}^K m(A_j)\right)^{1-r}}\geq C\frac{\sum_{i=1}^K a_i^r 
	m(A_i)}{(1\vee m([0,\infty)\times E))^{1-r}}, 
\]
which is \eqref{poiest2}.

For a general $(\calf_t)_{t\geq0}$-predictable process $H$,  one can choose a sequence $H_n$ of processes of the form 
\eqref{simple} such that $|H_n|\leq |H|$ for all $n\in\bbn$ and $H_n \to H$ as 
$n\to\infty$, pointwise in $(\om,t,x)$. If the right-hand side 
of \eqref{poiest} is finite, then inequality \eqref{poiest} follows from the 
dominated convergence theorem for stochastic integrals (see \cite[Equation~(2.6)]{Bichteler83}) on the left-hand side and for Lebesgue integrals on the right-hand side. If the right-hand side of \eqref{poiest} is infinite, then the estimates we have established for simple integrands, together with the BDG inequalities, imply that also the left-hand side of \eqref{poiest} is infinite.
\epr

\blem\label{split-pmoment} 
Suppose that $a\in\bbr$ and $X$ is a random variable with zero mean. Then for every $p\in(1,3]$, we have
\[ 
\bbe[|a+X|^p] \geq C_p(|a|^p + \bbe[|X|^p]) 
\]
where $C_p=1/4$ for $p \in (1,2]$ and $C_p=1/6$ for $p \in (2,3]$.
\elem

\begin{proof}
	First we prove the statement for $a =1$ and $p\in(1,2]$.
	The proof follows from the following simple inequalities:
	\begin{equation} \label{eq:lemma-auxineq}
		\begin{split}
			& (y-1)^p \geq \frac{1}{3} \left( y^p - 2 y + 1 \right), \quad y \geq 1, \\
			& (1-y)^p \geq \frac{1}{3} \left( y^p - 2 y + \frac{3}{4} \right), \quad y \in [0,1], \\
			& (y+1)^p \geq \frac{1}{3} \left( y^p + 2 y + 1 \right), \quad y \geq 0.
		\end{split}
	\end{equation}
	Indeed, denoting the distribution function of $X$ by $F$, (\ref{eq:lemma-auxineq}) and $\E [X] = 0$ imply
	\[
	\begin{split}
	\E [|1 + X|^p] & = \int_{-\infty}^\infty |1 + y|^p\, F(\dd y) \\
	& \geq  \int_{-\infty}^{-1} \frac{1}{3} \left( (-y)^p - 2 (-y) + 1 \right)\, F(\dd y) +
	\int_{-1}^0 \frac{1}{3} \left( (-y)^p - 2 (-y) + \frac{3}{4} \right)\, F(\dd y) \\
	& \quad+ \int_0^\infty \frac{1}{3} \left( y^p + 2 y + 1 \right)\, F(\dd y) \\
	& \geq \frac{1}{3} \E [|X|^p] + \frac23 \E [X] + \frac{1}{4} \geq \frac{1}{4} \left(  \E [|X|^p] + 1 \right).
	\end{split}
	\]
	For general $a\in\R$, the statement follows from
	\[
	\begin{split}
	\E [| a + X |^p] & = |a|^p \, \E \left[\left| 1 + \frac{X}{a}\right|^p\right]  \geq   \frac{|a|^p}{4}  \left( \frac{\E [|X|^p]}{|a|^p} + 1  \right)  = \frac{1}{4} \E [|X|^p] + \frac{1}{4} |a|^p.
	\end{split}
	\]
	
	Here is the proof of (\ref{eq:lemma-auxineq}). 
	The first inequality holds for $y=1$, and 
	\[
	\begin{split}
	p (y - 1)^{p-1} \geq p ( y^{p-1} -1) = \frac{1}{3} ( p y^{p-1} - 2) +
	\frac{2p}{3} y^{p-1} + \frac{2}{3} - p \geq \frac{1}{3} ( p y^{p-1} - 2),
	\end{split}
	\]
	that is the derivative of the left-hand side is greater than that of the right-hand side
	for all $ y \geq 1$. Thus the first inequality follows. For the second, using $y^p + 
	(1-y)^p \leq 1$, $y \in [0,1]$, we have
	\[
	\begin{split}
	3 (1-y)^p -  y^p + 2 y - \frac{3}{4} 
 &\geq 3 (1-y)^p - ( 1 - (1-y)^p) + 2 y - \frac{3}{4}  \geq 4 (1-y)^2 + 2 y - \frac{7}{4} \\
	& = \left( 2 y - \frac{3}{2} \right)^2,
	\end{split}
	\]
	which is nonnegative, so the second inequality is proved. Finally, for $y \geq 0$
	\[
	( y + 1 )^p \geq  y^p + 1 = \frac{1}{3} (y^p + 2 y + 1) +
	\frac{2}{3} (y^p  - y + 1) \geq \frac{1}{3} (y^p + 2 y + 1).
	\]
	
The proof is similar for $p\in (2,3]$, once the inequalities
\begin{equation} \label{eq:lemma-auxineq2}
\begin{split}
& (y-1)^p \geq \frac{1}{6} \left( y^p - 6 y + 1 \right), \quad y \geq 1, \\
& (1-y)^p \geq \frac{1}{3} \left( y^p - 3 y + 1 \right)
\geq \frac{1}{6} \left( y^p - 6 y + 1 \right), \quad y \in [0,1], \\
& (y+1)^p \geq \frac{1}{3} \left( y^p + 3 y + 1 \right)
\geq \frac{1}{6} \left( y^p + 6 y + 1 \right) , \quad y \geq 0,
\end{split}
\end{equation}
are established. We leave the proof of 
(\ref{eq:lemma-auxineq2}) to the interested reader.
\end{proof}

\bpr[\bff{\emph{Proof of Theorem~\ref{interm}}}] \benu
\item We assume $d\geq2$ here as the case $d=1$ will be treated in part (4). In particular, $p$ is always less than $2$ and $\La$ contains no Gaussian part. By 
Lemma~\ref{split-pmoment} and the BDG inequalities, we have for all 
$p\in(1,1+2/d)$
\beq\label{est}
\E[|Y(t,x)|^p] \geq C_p\left(L_f^p+ \E\left[\left(\iiint_0^t g^2(t-s,x-y)\si^2(Y(s,y)) z^2 
\,\mu(\dd s,\dd y,\dd z) \right)^{\frac{p}{2}}  \right]\right).\eeq
This estimate remains valid if we replace $\mu$ on the right-hand side by the measure 
\beq\label{mu-choice}
\mu_{\eps,\delta}^{(t,x)} (\dd s,\dd y,\dd z)  := 
\bone_{[0,t]}(s)\bone_{\{ g(t-s,x-y)>\eps \}}
\bone_{[-\delta,\delta]^\Comp}(z) 
\,\mu(\dd s,\dd y,\dd z)
\eeq
where $\eps>0$ is arbitrary and $\delta>0$ is chosen small enough such that 
$\la([-\delta,\delta]^\Comp)>0$. The corresponding intensity measure is given by 
\[
\nu_{\eps,\delta}^{(t,x)}(\dd s,\dd y,\dd z)  := 
\bone_{[0,t]} (s) \bone_{\{ g(t-s,x-y)>\eps \}}
\bone_{[-\delta,\delta]^\Comp}(z) \, \dd s\,\dd y\,\la(\dd z),
\]
and satisfies 
\begin{align*}
\nu^{(t,x)}_{\eps,\delta}([0,\infty)\times\R^d\times\R) &= \la([-\delta,\delta]^\Comp) \iint_0^t  
\bone_{\{g(s,y)>\eps\}} \,\dd s\,\dd y\\ 
&\leq \la([-\delta,\delta]^\Comp) \iint_0^\infty  
\bone_{\{g(s,y)>\eps\}} \,\dd s\,\dd y <\infty,
\end{align*}
with an upper bound independent of $(t,x)$. 
By Lemma~\ref{Poi-ineq}, we obtain (keeping in mind that $L_f>0$ and 
$L_\si>0$, and using the BDG inequality from the first to the second line)
\beq \label{lower-p}
\begin{split} 
	\E[|Y(t,x)|^p] 
	&\geq C_p \left( 1+ \E \left[ \left(\iiint_0^t g^2(t-s,x-y)\si^2(Y(s,y)) z^2 
	\,\mu^{(t,x)}_{\eps,\delta}(\dd s,\dd y,\dd z) \right)^{\frac{p}{2}}  \right]\right) \\
	& \geq C_p \left( 1 + \E \left[ \left| \iiint_0^t g(t-s,x-y)\si(Y(s,y)) z 
	\,( \mu^{(t,x)}_{\eps,\delta} - \nu^{(t,x)}_{\eps,\delta}) (\dd s,\dd y,\dd z) \right|^{p} \right] 
	\right)\\
	& \geq C_p \Bigg( 1 + \frac{\int_\R |z|^p\bone_{\{|z|>\delta\}}\,\la(\dd z)} 
	{(1\vee\la([-\delta,\delta]^\Comp)
		\iint_0^\infty  \bone_{\{g(s,y)>\eps\}} \,\dd s\,\dd y)^{1-\frac{p}{2}}} \\
	& \quad \times \iint_0^t g^p(t-s,x-y)\bone_{\{ g(t-s,x-y)>\eps \}} 
	\bbe[|Y(s,y)|^p] \,\dd s\,\dd y \Bigg)
\end{split}
\eeq
with a constant $C_p$ independent of $(t,x)$. As a consequence, the function 
\[
I_p(t):=\inf_{x\in\bbr^d} \bbe[|Y(t,x)|^p]
\]
satisfies 
\beq \label{renewal} 
I_p(t)\geq a_p + \int_0^t w_p(t-s)I_p(s)\,\dd s 
\eeq
for some $a_p>0$ where 
\beq\label{wdef}
w_p(t) = C_p \frac{\int_\R |z|^p\bone_{\{|z|>\delta\}} \,\la(\dd z)} 
{ (1\vee\la([-\delta,\delta]^\Comp) 
	\iint_0^\infty  \bone_{\{g(t,x)>\eps\}} \,\dd t\,\dd x)^{1-\frac{p}{2}}} 
\int_{\R^d} g^p(t,x) \bone_{\{ g(t,x)>\eps \}} \,\dd x.
\eeq
Recall from Lemmata~\ref{Poi-ineq} and \ref{split-pmoment} that both $C_p$ and $a_p$ 
can be assumed to be bounded away from $0$ if $p$ is bounded away from $1$. Since $g\notin L^{1+2/d}([0,T]\times\R^d)$ for any $T>0$ (cf.\ the calculations before Theorem~\ref{nomoments}), and the heat kernel decays exponentially in space, we have 
$\iint_0^\infty g^p(t,x)\bone_{\{ g(t,x)>\eps \}}\,\dd t\,\dd x \to \infty$ as 
$p\to 1+2/d$ and consequently,
\beq\label{explode}
\lim_{p\to 1+\frac{2}{d}} \int_0^\infty w_p(t)\,\dd t = \infty.
\eeq
Hence, there exists $p_0\in(1,1+2/d)$ such that $\int_0^\infty w_{p_0}(t)\,\dd 
t>1$. By classical renewal theory, see e.g.\  \cite[Theorem~V.7.1]{Asmussen03}, it follows that the solution to the equation 
\[ i(t)= a_p + \int_0^t w_{p_0}(t-s)i(s)\,\dd s  \] 
satisfies $i(t)\geq e^{\ga t}$ for all $t\geq t_0$ and some $t_0>0$ and $\ga>0$. Since $I_{p_0}(t)\geq i(t)$ by
 \cite[Theorem~7.11]{Khos}, we conclude that $\un\ga(p_0)>0$, and by Jensen's inequality, 
also $\un\ga(p)>0$ for all $p_0\leq p<1+2/d$.
\item Again, we only consider the case $d\geq2$.  A direct computation shows that 
\beq\label{pintegral}
\begin{split}
	\iint_0^\infty g^p(t,x)\bone_{\{ g(t,x)>\eps \}}\,\dd t\,\dd x &=\frac{2\pi^{\frac{d}{2}}}{\Ga(\frac{d}{2})} 
	\int_0^{\frac{1}{2\pi\kappa\eps^{2/d}}} \int_0^{\sqrt{-d\kappa t\log(2\pi\kappa\eps^{2/d} t)}} 
	\frac{e^{-\frac{pr^2}{2\kappa t}}}{(2\pi \kappa t)^{\frac{pd}{2}}} r^{d-1}\,\dd r\,\dd t\\
	&=\frac{\pi^{\frac{d}{2}}}{\pi\Ga(\frac{d}{2})\kappa\eps^{\frac{2}{d}-p}} 
\int_0^1\int_0^{\sqrt{\frac{-ds\log(s)}{2\pi \eps^{2/d}}}}   
s^{-\frac{pd}{2}}{e^{-\frac{pr^2\pi\eps^{2/d}}{s}}} r^{d-1}\,\dd r\,\dd s \\
	&=\frac{(\frac{d}{2})^{\frac{d}{2}}}{\pi\Ga(\frac{d}{2})\kappa\eps^{1+\frac{2}{d}-p}} \int_0^1\int_0^1 s^{\frac{pd(z^2-1)}{2}}(-s\log(s))^{\frac{d}{2}}z^{d-1}\,\dd z\,\dd s\\
	&= \frac{C_p}{\kappa\eps^{1+\frac{2}{d}-p}}
\end{split}
\eeq
for all $p\in(0,1+2/d)$.
This formula is still valid for $p=0$. Thus, for the function in \eqref{wdef}, which we 
denote by $w_\kappa(t)$ now since $\kappa$ is the parameter that interests us, there 
exists $C>0$ that is independent of $\kappa$ such that 
\beq\label{kappa-order} \lim_{\kappa\to 0}\int_0^\infty w_\kappa(t)\,\dd t = C 
\lim_{\kappa\to 0} \frac{\kappa^{1-\frac{p}{2}}}{\kappa} = \lim_{\kappa\to 
0}C\kappa^{-\frac{p}{2}} = \infty. \eeq
The proof can now be completed as in the first part of the theorem.
\item This part follows as before because $L_\si$ enters $C_p$ in \eqref{wdef} in a multiplicative way.
\item In dimension $1$ it suffices by Jensen's inequality to consider $p\in(1,2)$. Furthermore, by Lemma~\ref{split-pmoment} and the BDG inequalities, we may assume $\rho=0$ without loss of generality. Then the proof of part 
(1) remains valid up to equation \eqref{wdef}. Instead of varying the value of $p$, we 
now let $\eps\to 0$, keeping $p\in(1,2)$, $\kappa>0$ as well as $\delta>0$ fixed. Writing $w_\eps(t)$ in 
the following instead of $w_p(t)$ for the function in \eqref{wdef}, it follows from  \eqref{pintegral} that  
\[ 
\lim_{\eps\to0} \int_0^\infty w_\eps(t)\,\dd t = C \lim_{\eps\to0} 
\frac{\eps^{3(1-\frac{p}{2})}}{\eps^{3-p}}  = C\lim_{\eps\to0} \eps^{-\frac{p}{2}} = \infty,
\]
and the assertion follows.
\eenu
\epr

\bpr[\bff{\emph{Proof of Theorem~\ref{intfronts}}}] Let $\al>0$, $p\in(1,2\wedge (1+2/d))$ and write $x=(x_1,\ldots,x_d)$. By Proposition~\ref{thm:moments}, 
$x\mapsto \bbe[|Y(t,x)|^p]$ is integrable, so we deduce from  \eqref{lower-p} (with 
$\eps,\delta>0$ sufficiently small) and the hypothesis $L_\si>0$,
\begin{align*}
	\int_{x_1\geq \al t} \bbe[|Y(t,x)|^p]\,\dd x &\geq \frac14\Bigg(\int_{x_1\geq \al t} 
	|Y_0(t,x)|^p\,\dd x \\
	&\quad+ C \int_{x_1\geq \al t}\left(\iint_0^t 
	g^p(t-s,x-y)\bone_{\{g(t-s,x-y)>\eps\}} \bbe[|Y(s,y)|^p]\,\dd s\,\dd y \right)\,\dd 
	x\Bigg),
\end{align*}
where the constant 
$C$ is given by 
\beq
C = \int_\bbr |z|^p\bone_{\{|z|>\delta\}}\,\la(\dd z) \left( 1 \vee
\la([-\delta,\delta]^\Comp) \iint_0^\infty \bone_{\{g(t,x)>\eps\}}\,\dd t\,\dd x 
\right)^{\frac{p}{2}-1}.
\label{C-value}\eeq
Let us write 
$v(t):=\int_{x_1\geq \al t} \bbe[|Y(t,x)|^p]\,\dd x$, 
$v_0(t):=\int_{x_1\geq \al t} |Y_0(t,x)|^p\,\dd x$ and 
\[
h(t):=\int_{x_1\geq \al t} g^p(t,x)\bone_{\{g(t,x)>\eps\}}\,\dd x. 
\]
Using that $x_1-y_1\geq\al(t-s)$ and $y_1\geq \al s$ imply $x_1\geq \al t$, we obtain
\[
v(t)\geq \frac14\left(v_0(t) + C\int_0^t h(t-s)v(s)\,\dd s\right)
\]
for all $t\geq0$. A straightforward extension of \cite[Lemma~4.2]{Chen15} to the 
$d$-dimensional setting shows that $v_0(t)>0$ for all $t>0$. So on the one hand, if 
\beq\label{toachieve} 
C \int_0^\infty h(t)\,\dd t> 4, 
\eeq
it follows from renewal theory (see the proof of Theorem~\ref{interm}) that
\beq\label{limsupinf} 
\limsup_{t\to\infty} e^{-\beta t} v(t) =\limsup_{t\to\infty} e^{-\beta t}
\int_{x_1\geq \al t} \bbe[|Y(t,x)|^p]\,\dd x =\infty 
\eeq
whenever $\beta>0$ is sufficiently small. On the other hand, from 
Proposition~\ref{thm:moments}, we know that
\begin{align*} 
	\int_{x_1\geq \al^\prime t} \bbe[|Y(t,x)|^p]\,\dd x & \leq C e^{\beta^\prime t}\int_{x_1\geq \al^\prime t} e^{- c|x|} \,\dd x \leq C e^{\beta^\prime t}\int_{x_1\geq \al^\prime t} 
	e^{- \frac{c}{\sqrt{d}}(|x_1|+\ldots+|x_d|)} \,\dd x \\ &\leq C e^{\beta^\prime t} \int_{\al^\prime t}^\infty 
	e^{-\frac{c}{\sqrt{d}} x_1}\, \dd x_1  \leq C e^{(\beta^\prime -\al^\prime \frac{c}{\sqrt{d}}) t}
\end{align*}
for all $\al^\prime>0$, some $\beta^\prime>0$ and some $C>0$ that is independent of $t$. Thus, the last expression decays 
exponentially whenever $\al^\prime$ is large enough, so in this case, \eqref{limsupinf} 
implies
\[
\limsup_{t\to\infty} e^{-\beta t} \int_{\al t \leq x_1 < \al^\prime t} 
\bbe[|Y(t,x)|^p]\,\dd x = \infty, 
\]
from which the assertion follows because
\[
\sup_{|x|\geq \al t} \bbe[|Y(t,x)|^p] \geq  \sup_{x_1 \geq \al t} \bbe[|Y(t,x)|^p] \geq 
((\al^\prime -\al)t)^{-1} \int_{\al t \leq x_1 < \al^\prime t} \bbe[|Y(t,x)|^p]\,\dd x.
\]

So the only thing left to show is that we can achieve \eqref{toachieve} by proper choices 
of the parameters involved. Since the heat kernel is radially symmetric, we have 
$\int_0^\infty h(t)\,\dd t \geq \frac{1}{2d} \int_0^\infty \tilde 
h(t)\,\dd t$ where
\beq\label{h-tilde} \tilde h(t) = \int_{|x|\geq \tilde\al t} g^p(t,x)\bone_{\{ g(t,x)>\eps \}}\,\dd x  \eeq
and $\tilde\al = \sqrt{d}\al$.
Using polar coordinates and changing variables $s=  2 \pi \kappa \epsilon^{2/d}t$ and
$u= r^2 p \pi \epsilon^{2/d} / s$, we obtain
\begin{equation} \label{eq:htilde-lower1}
\begin{split}
\int_0^\infty \tilde h(t)\,\dd t 
&= \int_0^{\frac{1}{2\pi\kappa\epsilon^{2/d}}} 
\int_{\bbr^d} g^p(t,x) \bone_{\big\{\tilde \alpha t \leq |x|< \sqrt{-2\kappa t\log( 
		\epsilon(2\pi\kappa t)^{d/2} )} \big \}}\,\dd x\,\dd t\\
&= \frac{2 \pi^{\frac d 2}}{\Gamma(\frac d 2)}
\int_0^{\frac{1}{2\pi\kappa\epsilon^{2/d}}} 
\int_0^\infty \frac{e^{-\frac{pr^2}{2\kappa t}}}{(2\pi\kappa t)^{\frac{pd}{2}}} \bone_{\big\{\tilde \alpha t \leq r < \sqrt{-2\kappa t\log( 
		\epsilon(2\pi\kappa t)^{d/2} )} \big \}}r^{d-1}\,\dd r\,\dd t\\
& = 
\frac{ \epsilon^{-(1+\frac{2}{d} - p)}}{\Gamma(\frac d 2) p^{\frac d 2} (2 \pi\kappa)}
\int_0^1 \int_0^\infty s^{-\frac{d}{2}(p-1)} e^{-u} u^{\frac{d}{2}-1}
\bone_{\left\{\frac{\tilde \alpha^2}{\kappa^2 \epsilon^{2/d}} \frac{s p}{4\pi} \leq 
	u \leq \frac{pd}{2} \log s^{-1} \right\}} \dd u \,  \dd s.
\end{split}
\end{equation}
Note that the latter integral depends on the parameters $\tilde \alpha$, 
$\kappa$, and $\eps$ only through the ratio
\begin{equation} \label{eq:def_R}
R = R(\tilde \alpha, \kappa, \epsilon) = \frac{\tilde \alpha^2}{\kappa^2 
	\varepsilon^{2/d}}. 
\end{equation}
Since $s$ is increasing and $\log s^{-1}$ is decreasing, 
\[
R \frac{p}{4 \pi} s \leq \frac{pd}{2} \log s^{-1}, \quad s \in (0, s_0), 
\]
with $s_0 = \min \{ 2 \pi d / R, e^{-1} \}$. Therefore, the integral in 
\eqref{eq:htilde-lower1} is not less than
\[
\begin{split}
& \int_0^{s_0} \int_{\frac{R p}{4\pi} s}^{\frac{pd}{2} \log s^{-1}}
s^{-\frac{d}{2}(p-1)} e^{-u} u^{\frac{d}{2} -1} \dd u \, \dd s 
\geq 
\int_0^{\frac{Rp}{4 \pi} s_0}  e^{-u} u^{\frac{d}{2}-1}
\int_0^{\frac{4 \pi}{Rp}u} s^{-\frac{d}{2}(p-1)} \dd s \, \dd u \\
& = 
\left( \frac{4 \pi}{ R p} \right)^{1 - \frac{d}{2}(p-1)}
\frac{2}{d(1 + \frac{2}{d}-p)} \,
\gamma \!\left(1 + d \left(1-\frac{p}{2} \right), \frac{R p s_0 }{4 \pi} \right),
\end{split}
\]
where the second inequality follows from Fubini's theorem, and
\begin{equation} \label{incomplete-gamma}
\gamma(x,T) = \int_0^T t^{x-1} e^{-t} \dd t, \quad x > 0,\quad T \geq 0, 
\end{equation}
stands for the lower incomplete gamma function. Substituting back into 
\eqref{eq:htilde-lower1}, we obtain
\begin{equation} \label{eq:htilde-lower2}\begin{split}
\int_0^\infty \tilde h(t)\,\dd t
&\geq 
 \frac{2^{2-d(p-1)}\kappa^{-1} \epsilon^{-(1+\frac{2}{d} - p)} R^{\frac{d}{2}(p-1)-1}}
{d\Ga(\frac d2)\pi^{\frac d 2 (p-1)}p^{1+d(1-\frac p 2)}(1 + \frac{2}{d}-p)} 
\gamma\left(1 + d \left(1-\frac{p}{2} \right), \frac{R s_0}{4 \pi} \right)\\
&= \frac{2(2\kappa)^{1-d(p-1)} \tilde\al^{-2(1-\frac{d}{2}(p-1))}}
{d\Ga(\frac d2)\pi^{\frac d 2 (p-1)}p^{1+d(1-\frac p 2)}(1 + \frac{2}{d}-p)} 
\gamma\left(1 + d \left(1-\frac{p}{2} \right), \frac{R s_0}{4 \pi} \right),
\end{split}
\end{equation}
where  $R$ is given in 
\eqref{eq:def_R}.

Consequently, when $\tilde\al^2\eps^{-2/d}\geq 2\pi\kappa^2$, we have
\beq\label{htildeintegral}
\begin{split}
	\int_0^\infty \tilde h(t)\,\dd t &\geq 
	\frac{2(2\kappa)^{1-d(p-1)}  
		\tilde 
		\al^{-2({1-\frac{d}{2}(p-1)})}\ga(1+d(1-\frac{p}{2}),\frac16) } {d\Ga(\frac d2)\pi^{\frac d 2 (p-1)}p^{1+d(1-\frac p 2)}(1 + \frac{2}{d}-p)}.
\end{split}
\eeq

Part (1) of the theorem in dimension $d\geq2$ now follows from the observation that the right-hand side of 
\eqref{htildeintegral} tends to $\infty$ as $p\to 1+2/d$, for any given $\kappa,\al>0$ and small 
values of $\eps$ and $\delta$ (note that the constant $C$ in \eqref{toachieve} is 
bounded for $p$ in a neighborhood of $1+2/d$). 

For (2) choose $\alpha = \sqrt{2 \pi} \epsilon^{1/d} \kappa$, with $\epsilon$ being fixed 
and $\kappa \to 0$. The lower bound in \eqref{htildeintegral} is of order $\kappa^{-1}$, 
while the constant $C$ in \eqref{C-value} is of order $\kappa^{1-p/2}$ by \eqref{pintegral}. 
Thus the statement follows by choosing $\kappa$ sufficiently small.
Similar considerations, compare also with the proof of Theorem~\ref{interm}, also show 
part (3) of the theorem. 

For part (4), i.e., if we are in dimension $d=1$,  let us 
assume $\rho=0$ without loss of generality and choose $\tilde\al = 
\sqrt{2\pi\kappa^2}\eps$, with $\kappa$ and $p$ being fixed this time. Then, for any 
given $p$ and $\kappa$, the right-hand side of \eqref{htildeintegral}
is of order $\eps^{-2(1-(p-1)/2)}=\eps^{-(3-p)}$ under the hypotheses of the 
theorem, while the constant $C$ in \eqref{toachieve} is of order $\eps^{3(1-p/2)}$ by 
\eqref{C-value} and \eqref{pintegral}, so we can achieve \eqref{toachieve} by taking 
$\eps$ small enough.
\epr

\bpr[\bff{\emph{Proof of Theorem~\ref{nomoments}}}]
Let $T>t_0$, $p=1+2/d$ and assume the opposite of \eqref{infmom}. Then by Minkowski's integral inequality,
\begin{align*} \left\| b\iint_0^t g(t-s,x-y)\si(Y(s,y)) \, \dd s\,\dd y \right\|_p&\leq |b|\iint_0^t g(t-s,x-y)\|\si(Y(s,y))\|_p\,\dd s\,\dd y\leq  C|b|T
\end{align*}
for all $(t,x)\in[0,T]\times\R^d$. Similarly, if $d=1$, and we have by the BDG inequality, 
\begin{align*} &\left\| \rho\iint_0^t g(t-s,x-y)\si(Y(s,y)) \, W(\dd s,\dd y) \right\|_p
	\leq C_p |\rho|\left\|\iint_0^t g^2(t-s,x-y)\si^2(Y(s,y))\,\dd s\,\dd y\right\|_{\frac{p}{2}}^{\frac{1}{2}}\\
	& \qquad\leq C_p |\rho| \left( \iint_0^t g^2(t-s,x-y)\|\si(Y(s,y))\|_p^2\,\dd s\,\dd y \right)^\frac12 \leq C_p |\rho| T^\frac14.
\end{align*}
Therefore, we deduce, if the left-hand side of \eqref{infmom} was finite, then
\beq\label{nonfin} \bbe\left[\left| \iiint_0^t g(t-s,x-y)\si(Y(s,y))z\,(\mu-\nu)(\dd s,\dd y,\dd z)\right|^p \right]<\infty \eeq
as well. We now show that this cannot be true. Indeed, because $\si(Y_0(t_0,x_0))\neq0$, there exists $(t_1,x_1)\in (0,t_0]\times\R^d$ with $\bbp[\si(Y(t_1,x_1))\neq0]>0$ (otherwise, we would have $\si(Y(s,y))=0$ for all $(s,y)\in(0,t_0]\times\R^d$ and hence $Y(t,x)=Y_0(t,x)$ for all $(t,x)\in[0,t_0]\in\R^d$ by \eqref{SHE}, which together contradict \eqref{sigma-not-zero}). Therefore, we have $\bbe[|\si(Y(t_1,x_1))|]>0$, and because the solution to \eqref{SHE} is $L^r$-continuous on $(0,\infty)\times\R^d$ for all $r<1+2/d$ (see \cite[Theorem~4.7(2)]{Chong}), there exist $\eps,\delta>0$ such that 
\beq\label{epsdelta} 0\leq t_1-s<\delta,~|x_1-y|\leq \sqrt{\delta} \implies  \bbe[|\si(Y(s,y))|] >\eps. \eeq

As seen in the proof of Theorem~\ref{interm}, the left-hand side of \eqref{nonfin} is 
greater than or equal to a constant times the same expression with $\mu$ and $\nu$ 
replaced by $\tilde\mu$ and $\tilde\nu$, respectively, where 
\[ 
\tilde\mu(\dd s,\dd y,\dd z) = \bone_{(t_1-\delta,t_1]}(s)\bone_{\{|x_1-y|\leq 
\sqrt{t_1-s}\}} (y) \bone_{\R\setminus[-a,a]}(z)\,\mu(\dd s,\dd y,\dd z), \]
and $\tilde \nu$ is its compensator. The number $a>0$ is chosen such that $\la(\R\setminus[-a,a])>0$. Now if $d\geq2$ (and consequently $p\leq 2$), Lemma~\ref{Poi-ineq} gives the estimate
\begin{align*}
	&\bbe\left[\left| \iiint_0^{t_1} g(t_1-s,x_1-y)\si(Y(s,y))z\,(\mu-\nu)(\dd s,\dd y,\dd z)\right|^p \right]\\
	&\qquad \geq C\frac{\iiint_0^{t_1} g^p(t_1-s,x_1-y)\bbe[|\si(Y(s,y))|^p] |z|^p \,\tilde\nu(\dd s,\dd y,\dd z)}{(1\vee \tilde\nu([0,t]\times \bbr^d\times\R))^{1-\frac{p}{2}}}\\
	&\qquad=C\left( \int_{0}^{\delta} \int_{|y|\leq \sqrt{s}} \int_{\bbr\setminus[-a,a]} \,\dd s\,\dd y\,\la(\dd z) \right)^{\frac{p}{2}-1}\\
	&\qquad \quad\times \int_{t_1-\delta}^{t_1}\int_{\R^d}\int_\R g^p(t_1-s,x_1-y)\bone_{\{|x_1-y|\leq \sqrt{t_1-s}\}}\bbe[|\si(Y(s,y))|^p] |z|^p\bone_{\{|z|>a\}}\,\dd s\,\dd y\,\la(\dd z)\\
	&\qquad\geq C\eps^p \delta^{(1+\frac{d}{2})(\frac{p}{2}-1)} \int_0^\delta \int_{|x|\leq \sqrt{t}} g^p(t,x)\,\dd t\,\dd x\\ 
	&\qquad= C\int_0^\delta \int_{|x|\leq \sqrt{t}} g^p(t,x)\,\dd t\,\dd x .
\end{align*}
The last line is a valid lower bound also in the case $d=1$, possibly with another value of $C$, as a consequence of \cite[Theorem~1]{MR}. But for $p=1+2/d$, we have
\[ 
\begin{split}
\int_0^\delta \int_{|x|\leq \sqrt{t}} g^p(t,x)\,\dd t\,\dd x
& = \int_0^\delta 
\int_{|x|\leq\sqrt{t}} \frac{e^{-\frac{p|x|^2}{2\kappa t}}}{(2\pi\kappa t)^{\frac{pd}{2}}} 
\,\dd t\,\dd x  \\
& \geq \frac{Ce^{-\frac{p}{2\kappa}}}{(2\pi\kappa )^{\frac{pd}{2}}} 
\int_0^\delta t^{-\frac{pd}{2}}t^{\frac{d}{2}} \,\dd t = C\int_0^\delta \frac1{t}\,\dd t= 
+\infty,
\end{split}
\]
proving that \eqref{nonfin} is wrong.
\epr

\bpr[\bff{\emph{Proof of Theorem~\ref{pos-mean}}}] For every $(t,x)\in(0,\infty)\times\bbr^d$, we have
\begin{align*}
\bbe[|Y(t,x)|]\geq \bbe[Y(t,x)] &= Y_0(t,x) + b\int_0^t\int_{\bbr^d} 
g(t-s,x-y)\bbe[\si(Y(s,y))]\,\dd s\,\dd y\\
& \geq Y_0(t,x) + 
b L_\si\int_0^t\int_{\bbr^d} g(t-s,x-y)\bbe[|Y(s,y)|]\,\dd s\,\dd y.
\end{align*}
Since the integral of $g$ on $(0,\infty)\times\R^d$ is infinite, the theorem follows from 
the renewal methods as used in Theorems~\ref{interm} and \ref{intfronts}.
\epr

\bpr[\bff{\emph{Proof of Theorem~\ref{interm-neg}}}] By Proposition~\ref{cable}, we can equally consider the stochastic cable equation \eqref{SHE-b} driven by the noise $\dot M$. 
\benu
\item We only carry out the proof for $c=0$; the arguments are similar for $c>0$ and we leave the details to the reader. Starting with $d\geq2$ and $p\in(1,1+2/d)$, by virtually the same calculations as in the proof of Theorem~\ref{interm}, the function $I_p(t)=\inf_{x\in\R^d} \E[|Y(t,x)|^p]$ satisfies \eqref{renewal} with $a_p$ replaced by the function $a_p(t)=a_p^\prime e^{-p|b|\si_0 t}$ where $a_p^\prime>0$ is a constant, and with $g$ replaced by $g^\prime$ in the definition \eqref{wdef} of $w_p$. But still, we have \eqref{explode}, so the renewal methods go through and the conclusion of Theorem~\ref{interm}(1) is valid. Statement (3) of the same theorem can be derived in a similar way. 

For statement (2), we observe that the truncated jump measure in \eqref{mu-choice} does not need to use the same kernel function as in \eqref{est} a priori (it only needs to have a finite intensity measure). Hence, for imitating the proof of Theorem~\ref{interm}(2), we only replace $g$ by $g'$ in \eqref{est}. For the indicator function $\bone_{\{g(t-s,x-y)>\eps\}}$ in \eqref{mu-choice}, by contrast, we replace $g(t,x)$ by $g(1;t,x)$ (i.e., the heat kernel with $\kappa=1$ and without the $e^{-|b|\si_0 t}$ factor). As a consequence, the function in \eqref{wdef} becomes
\[ w'_\kappa (t) = C_p \frac{\int_\R |z|^p\bone_{\{|z|>\delta\}} \,\la(\dd z)} 
{ (1\vee\la([-\delta,\delta]^\Comp) 
	\iint_0^\infty  \bone_{\{g(1; t,x)>\eps\}} \,\dd 
	t\,\dd x)^{1-\frac{p}{2}}} 
\int_{\R^d} e^{-p|b| \sigma_0 t} g^p (\kappa; t,x) \bone_{\{ g(1; t,x)>\eps \}} \,\dd x.  \]
Only the integral term in the previous line depends on $\kappa$. Hence, 
we conclude from \eqref{calculation} (the calculation there is valid up to the third line for any value of $\beta$) that 
$\int_0^\infty w_\kappa'(t)\,\dd t$ is of order $\kappa^{-(p-1)d/2}$.

For $d=1$, we need to let $p\to3$. The BDG inequalities allow us to ignore the Gaussian part, so by \cite[Theorem~1]{MR} and Lemma~\ref{split-pmoment}, we have that
\beq\label{d1}
\begin{split}
	\bbe[|Y(t,x)|^p] &\geq \frac16\Bigg( Y_0(t,x)^p + \si_0^p \bbe\left[ \left( \iiint_0^t | {g^\prime}(t-s,x-y)Y(s,y)z|^2\,\nu(\dd s,\dd y,\dd z)\right)^\frac{p}{2} \right] \\
	&\quad + \si_0^p \bbe\left[ \iiint_0^t | {g^\prime}(t-s,x-y)Y(s,y)z|^p\,\nu(\dd s,\dd y,\dd z) \right] \Bigg)\\
	&\geq  \frac16\left( Y_0(t,x)^p +  \si_0^p m^p_\la(p) \iiint_0^t  {g^\prime}^p(t-s,x-y)\bbe[|Y(s,y)|^p]\,\dd s\,\dd y\right),
\end{split}
\eeq
so we can complete the proof as in the case $d\geq2$ above.
\item Since $\ov\ga(p)<0$ implies $\ov\la(p)=0$, we can assume $c=0$ in this part of the theorem. Furthermore, by the hypotheses of Proposition~\ref{thm:moments}, the assumption that $m_\la(1+2/d)<\infty$, and Jensen's inequality, we may assume that $p$ is large enough such that $m_\la(p)$ is finite, and if $d=1$ and $\rho\neq0$, that $p\geq2$. Writing $C_{\beta, c}(b,\rho, \lambda, \kappa, p)$ for the constant 
$C_{\beta, c}(\kappa,p)$ in Proposition~\ref{prop:stochYoung-d} to stress the dependence of the constant on the 
other parameters, we obtain with identical calculations as in the proof of Proposition~\ref{prop:stochYoung-d} that 
$\| g' \circledast \Phi \|_{p, \beta, 0} \leq 
C_{\beta + |b| \sigma_0, 0}(0, \rho, \lambda, \kappa, p)$.
In particular, if we re-examine the proof of Proposition~\ref{thm:moments} and the formula \eqref{Cbetagamma}, we see that whenever $\kappa$ or $|b|$ is large, or $\si_0$ is small, there exists $\beta<0$ such that $C_{\beta+|b|\si_0,0}(0,\rho,\lambda,\kappa,p)<1/\si_0$ and thus $\|Y\|_{p,\beta,0}<\infty$ and $\ov\ga(p)<0$.
\eenu
\epr

\bpr[\bff{\emph{Proof of Theorem~\ref{compprin}}}] Let us introduce the following truncations of $\La$:
\begin{align*} \La_n(\dd t,\dd x) &= b\psi_n(x)\,\dd t\,\dd x +\psi_n(x)\int_0^\infty 
z\bone_{\{z>\frac{1}{n} \}}\,(\mu-\nu)(\dd t,\dd x,\dd z)\\
	&= \left(b -  \int_0^\infty z \bone_{\{z>\frac{1}{n}\}}\,\la(\dd z)\right) \psi_n(x)\,\dd t\,\dd x+ \psi_n(x)\int_0^\infty z \bone_{\{z>\frac{1}{n}\}}\,\mu(\dd t,\dd x,\dd z)\\
	&=: b_n \psi_n(x)\,\dd t\,\dd x + \La^+_n(\dd t,\dd x),\quad n\in\N,  \end{align*}
where $\psi_n(x)=\psi(|x|/n)$ and $\psi\colon[0,\infty)\to[0,1]$ is a smooth function with $\bone_{[0,1]}\leq \psi\leq \bone_{[0,2]}$.
If $Y_n$ denotes the solution to \eqref{SHE} with noise $\La_n$, we have by \cite[Theorem~1]{Chen16} that $Y_n(t,x)\to Y(t,x)$ in $L^p$ for all $(t,x)\in[0,\infty)\times\R^d$ (the cited result remains valid for the smooth truncation functions $\psi_n$ instead of the indicator functions $\bone_{[-n,n]^d}$). So if we can show that almost surely, with obvious notation, $Y_{n,1}(t,x)\geq Y_{n,2}(t,x)$  for all $(t,x)\in[0,\infty)\times\R^d$, then it follows that $Y_1(t,x)\geq Y_2(t,x)$ for all $(t,x)\in[0,\infty)\times\R^d$ upon choosing separable modifications of $Y_1$ and $Y_2$, which is always possible, see \cite[Theorem~II.2.4]{Doob53}. 

Now notice that for every $T>0$, the measure $\La_n^+$  only has finitely many jumps on $[0,T]\times\R^d$ almost surely. Let $T_0=0$ and $(T_i,X_i,Z_i)$, $i=1,\ldots,N_n(T)$, be the corresponding jump times, positions and sizes. The crucial observation is now that between $(T_{i-1},T_i)$, in absence of jumps, both $Y_{n,1}$ and $Y_{n,2}$ satisfy the deterministic PDE
\[ \partial_t Y_{n,j} (t,x)=\frac\kappa 2 \Delta Y_{n,j}(t,x)+b_n\si(Y_{n,j}(t,x))\psi_n(x),\quad j=1,2, \]
respectively.
Since $f_1\geq f_2$ and $\si$ is Lipschitz continuous, the comparison principle for the deterministic heat equation (see  \cite[Theorem~II]{Besala63}) implies  $Y_{n,1}(t,x)\geq Y_{n,2}(t,x)$ for all $(t,x)\in[0,T_1)\times\R^d$.
By induction, we may therefore assume that $Y_{n,1}(t,x)\geq Y_{n,2}(t,x)$ holds for all $(t,x)\in [0,T_{i})\times\R^d$ and then prove the same relation for $(t,x)\in[T_{i},T_{i+1})\times\R^d$. But since $Z_i\geq0$ and hence
\begin{align*} Y_{n,1}(T_i,x)&=Y_{n,1}(T_i-,x)+\si(Y_{n,1}(T_i-,X_i))Z_i\delta_{X_i}(x)\\
&\geq Y_{n,2}(T_i-,x)+\si(Y_{n,2}(T_i-,X_i))Z_i\delta_{X_i}(x)=Y_{n,2}(T_i,x) \end{align*}
by the induction hypothesis and the monotonicity property of $\si$, this again follows from the deterministic comparison principle (by considering smooth approximations of the Dirac delta function, the result of \cite{Besala63} extends to the measure-valued initial conditions encountered here).

Concerning the second statement of the theorem, the nonnegativity of $Y$ follows from the first part by comparison with the zero solution corresponding to a zero initial condition. Next, observe that the mean function $m(t,x)=\bbe[Y(t,x)]$ satisfies $m(0,x)=f(x)$ and 
\beq\label{mean-pde}  \partial_t m(t,x)=\Delta m(t,x) +b\bbe[\si(Y(t,x))]\begin{cases} \leq \Delta m(t,x) + (b\vee 0) L m(t,x),\\ \geq \Delta m(t,x) + (b\wedge 0) L m(t,x).\end{cases} \eeq
Again by the deterministic comparison principle, $m'(t,x)\leq m(t,x)\leq m''(t,x)$ where $m'$ (resp.\ $m''$) is the solution to \eqref{mean-pde} with equality instead of ``$\geq$'' (resp.\ ``$\leq$''). Since $m'$ (resp.\ $m''$) is given by the left-hand side (resp.\ right-hand side) of \eqref{mean-calc}, all assertions follow.
\epr

\subsection{Proofs for Section~\ref{Sect4}}

\bpr[\emph{\bff{Proof of Theorem~\ref{asymptotics-1}}}]  \benu\item If $\beta_0>0$ satisfies \eqref{satisfies}, then the proof of Proposition~\ref{thm:moments} reveals that $\|Y\|_{p,\beta_0,0}<\infty$, and hence $\ov\ga(p) \leq p\beta_0$.
When $\la\not\equiv0$ and $p$ is close enough to $1+2/d$, the second summand in \eqref{Cbetagamma} is always the term of leading order. Thus, \eqref{satisfies} holds as soon as $\beta_0$ satisfies
\[ 
\frac{\Ga(1-(p-1)\frac{d}{2})^{\frac{1}{p}}}{\beta_0^{\frac{1}{p}-\frac{d}{2p}(p-1)}} < C 
\iff \beta_0> C^{-\frac{2/d}{ 1+2/d-p }} 
\textstyle\Ga \left(\frac{d}{2}\left( 1+\frac{2}{d}-p \right) 
\right)^{\frac{2/d}{ 1+2/d-p }}
\]
for some finite constant $C$ independent of $p$. Since $x\Ga(x) =\Ga(1+x)\to 1$ as 
$x\to0$, we can choose 
\[ 
\beta_0=C^{-\frac{2/d}{ 1+2/d-p }}\left( 
\frac{\frac{2}{d}}{ 1+\frac{2}{d}-p 
	}\right)^{\frac{2/d}{ 1+2/d-p}}
\]
when $p$ is sufficiently close to $1+2/d$, which implies 
\beq\label{asymp-eq1} \limsup_{p\to 1+\frac{2}{d} }\frac{1+\frac{2}{d}-p }{\left|\log\left(1+\frac{2}{d}-p  \right)\right|}
\log \ov\ga(p)\leq \frac{2}{d}.\eeq
The upper bound in \eqref{kappato0} follows similarly. 

For the lower bounds in \eqref{p-upperbound} and \eqref{kappato0}, we first consider the case $b=0$. 
For $d\geq2$ let $\beta_1=\beta_1(p)$ be the number for 
 which
 \[
 \int_0^\infty w_p(t)e^{-\beta_1 t}\,\dd t = 1
 \]
 where $w_p$ is given by \eqref{wdef}. Recalling \eqref{incomplete-gamma}, and assuming that $p$ is close to $1+2/d$, and $\eps,\delta >0$ are small enough such that \eqref{lower-p} holds, we have that 
 \beq\label{beta-integral}
 \begin{split}
 	&\iint_0^\infty e^{-\beta t}g^p(t,x)\bone_{\{ g(t,x)>\eps \}}\,\dd t\,\dd x \\
 	&\qquad= \int_0^{\frac{1}{2\pi\kappa \eps^{2/d}}} \frac{e^{-\beta t}}{(2\pi\kappa t)^{\frac{pd}{2}}} 
 	\int_{\bbr^d} e^{-\frac{p|x|^2}{2\kappa t}}\bone_{\{ |x|^2<-2\kappa t\log( \eps(2\pi\kappa t)^{d/2} ) \}} 
 	\,\dd x\,\dd t\\
 	&\qquad=\frac{2\pi^{\frac{d}{2}}}{\Gamma(\frac{d}{2})} \int_0^{\frac{1}{2\pi \kappa\eps^{2/d}}} 
 	\frac{e^{-\beta t}}{(2\pi\kappa t)^{\frac{pd}{2}}} \int_0^{\sqrt{-2\kappa t\log( \eps(2\pi\kappa t)^{d/2} )}} 
 	e^{-\frac{pr^2}{2\kappa t}} r^{d-1}\,\dd r\,\dd t\\
 	&\qquad=\frac{1}{p^\frac{d}{2}\Ga(\frac{d}{2})(2\pi\kappa)^{\frac{d}{2}(p-1)}} 
 	\int_0^{\frac{1}{2\pi\kappa \eps^{2/d}}} \frac{e^{-\beta t}}{t^{\frac{d}{2}(p-1)}} \ga\left( 
 	\textstyle\frac{d}{2},-p\log\left( \eps(2\pi\kappa t)^{\frac{d}{2}} \right) \right) \,\dd t\\
 	&\qquad \geq 
 	\frac{\ga(\frac{d}{2},1)}{p^{\frac{d}{2}}\Ga(\frac{d}{2})(2\pi\kappa)^{\frac{d}{2}(p-1)}} 
 	\int_0^{\frac{e^{-2/(pd)}}{2\pi\kappa \eps^{2/d}}} \frac{e^{-\beta t}}{t^{\frac{d}{2}(p-1)}}  
 	\,\dd t\\
 	&\qquad = \frac{\ga(\frac{d}{2},1) \ga(1-\frac{d}{2}(p-1), (2\pi\kappa)^{-1} 
 		\eps^{-\frac{2}{d}}e^{-\frac{2}{pd}}\beta)}{p^{\frac{d}{2}}\Ga(\frac{d}{2})(2\pi\kappa)^{\frac{d
 			}{2}(p-1)}\beta^{1-\frac{d}{2}(p-1)}}.
 \end{split}\eeq
 It follows for $\beta\geq 2\pi \kappa e^{2/(pd)} \eps^{2/d}$ 
 that
 \begin{align}\label{beta-lowerbound}
 	\begin{split} \iint_0^\infty e^{-\beta t}g^p(t,x)\bone_{\{ g(t,x)>\eps \}}\,\dd t\,\dd x 
 		&\geq \frac{\ga(\frac{d}{2},1) \ga(1-\frac{d}{2}(p-1), 
 			1)}{p^{\frac{d}{2}}\Ga(\frac{d}{2})(2\pi\kappa)^{ \frac{d}{2}(p-1) 
 			}\beta^{1-\frac{d}{2}(p-1)}} \\
 		& \geq \frac{\ga(\frac{d}{2},1)(1-e^{-1})} 
 		{p^{\frac{d}{2}} \Ga(\frac{d}{2}) (2\pi\kappa)^{\frac{d}{2}(p-1)} 
 			(1-\frac{d}{2}(p-1)) \beta^{1-\frac{d}{2} (p-1)}}, 
 	\end{split}
 \end{align}
where the last step uses $\ga(1,1) = 1-e^{-1}$ and the fact that $x\ga(x,1)$ is a continuous decreasing function on 
$[0,1]$. Indeed, the latter follows from the identity $x\ga(x,1)=\ga(x+1,1)+e^{-1}$, which can be proved by integration by parts. Observing that the factor in front of the integral 
in \eqref{wdef} is bounded for $p$ around $1+2/d$, we deduce from 
 \eqref{beta-lowerbound} that
 \beq\label{xtox} 
 \beta_1 \geq \left( \frac{C}{1-\frac{d}{2}(p-1)} \right)^{\frac{1}{1-d(p-1)/2}} = 
 \left( C\frac{\frac{2}{d}}{1+\frac{2}{d}-p} 
 \right)^{\frac{2/d}{1+2/d-p}} 
 \eeq
 for some  constant $C$  independent of $p$.
 Hence we obtain from \cite[Theorem~V.7.1]{Asmussen03} that
 \[
 \un\ga(p) \geq 
 \beta_1 \geq \left( C\frac{\frac{2}{d}}{1+\frac{2}{d}-p} 
 \right)^{\frac{2/d}{1+2/d-p}}, 
 \]
 which implies  
 \beq\label{asymp-eq2}  \liminf_{p\to 1+\frac{2}{d} }\frac{1+\frac{2}{d}-p }{\left|\log\left(1+\frac{2}{d}-p  \right)\right|}
 \log \un\ga(p)\geq \frac{2}{d} \eeq
 and hence \eqref{p-upperbound} together with \eqref{asymp-eq1}. For $d=1$, if we estimate as in \eqref{d1}, the same arguments apply and only some constants would change that have no impact on the result. 
 
 For the lower bound in \eqref{kappato0}, the estimates \eqref{beta-integral} and \eqref{beta-lowerbound} can be re-used in principle, but we need to make a small change in our arguments because the denominator in \eqref{wdef} involves the kernel $g$ and therefore the parameter $\kappa$, which would lead to a suboptimal lower bound. In order to avoid this, we proceed as in the proof of Theorem~\ref{interm-neg}(2), and construct the measure in \eqref{mu-choice} by using the indicator function $\bone_{\{ g(1; t-s,x-y)>\eps \}}$ instead of $\bone_{\{ g(t-s,x-y)>\eps \}}$, where $g(1;t,x)$ is the heat kernel with $\kappa=1$. Then we have for $\kappa\leq1$ and $\beta\geq 2\pi\eps^{2/d}e^{2/(pd)}$,
 \beq\begin{split}
 	&\iint_0^\infty e^{-\beta t} g^p(t,x)\bone_{\{g(1;t,x)>\eps\}}\,\dd t\,\dd x\\ 
 	&\qquad=\frac{1}{p^\frac{d}{2}\Ga(\frac{d}{2})(2\pi\kappa)^{\frac{d}{2}(p-1)}} 
 	\int_0^{\frac{1}{2\pi \eps^{2/d}}} \frac{e^{-\beta t}}{t^{\frac{d}{2}(p-1)}} \ga\left( 
 	\textstyle\frac{d}{2},-p\kappa^{-1}\log\left( \eps(2\pi t)^{\frac{d}{2}} \right) \right) \,\dd t\\
 	&\qquad \geq 
 	\frac{\ga(\frac{d}{2},1)}{p^{\frac{d}{2}}\Ga(\frac{d}{2})(2\pi\kappa)^{\frac{d}{2}(p-1)}} 
 	\int_0^{\frac{e^{-2/(pd)}}{2\pi \eps^{2/d}}} \frac{e^{-\beta t}}{t^{\frac{d}{2}(p-1)}}  
 	\,\dd t
  \geq \frac{\ga(\frac{d}{2},1) \ga(1-\frac{d}{2}(p-1), 1)}{p^{\frac{d}{2}}\Ga(\frac{d}{2})(2\pi\kappa)^{\frac{d}{2}(p-1)}\beta^{1-\frac{d}{2}(p-1)}}.
  \end{split} \label{calculation}
\eeq
Thus, $\beta \geq C \kappa^{-\frac{p-1}{1+2/d- p}}$, proving the lower bound in 
\eqref{kappato0}.

Now let us explain why the proof of the lower bounds, for both $p\to1+2/d$ and $\kappa\to0$, remains essentially unchanged for $b<0$ or $b>0$. Indeed, if $\si$ is given by \eqref{PAM}, Proposition~\ref{cable} implies that we have to multiply $g$ by a factor $e^{b\si_0 t}$. But under the truncation $\bone_{\{g(t,x)>\eps\}}$ (resp. $\bone_{\{g(1;t,x)>\eps\}}$ when $\kappa\to0$ is considered), we have $t<T$ where $T=(2\pi\eps^{2/d})^{-1}$ is independent of $p$ (resp. $\kappa$). In particular, $g$ and $ge^{b\si_0 t}$ differ at most by a multiplicative constant $e^{b\si_0 T}$ on $[0,T]$, which is irrelevant for the calculations above. 

\item The upper bound for $\ov\la(p)$ in \eqref{p-upperbound-2} as $p\to 1+2/d$ 
follows from (1) because we have \eqref{help}. For the upper bound in 
\eqref{kappato0-2}, observe from \eqref{help} that $\ov\la(p) \leq \beta_0/c$ where 
$\beta_0$ was introduced in the proof of Proposition~\ref{thm:moments}. Upon inspection of 
formula \eqref{Cbetagamma}, we see that $\beta_0$ must satisfy
\[ \frac{C}{\beta_0-\frac12 \kappa c^2 d}+ 
\frac{C^\prime}{\kappa^{\frac{d(p-1)}{2p}}(\beta_0-\frac12 \kappa c^2 
d)^\frac{2-d(p-1)}{2p}} + \frac{C^{\prime\prime}}{\kappa^{\frac14}(\beta_0-\frac12\kappa 
c^2)^{\frac14}}\bone_{\{d=1,\ p\geq2\}} \leq 1. \]
As long as $\la\not\equiv0$, the second summand is the dominant one for small $\kappa$, so 
$\beta_0$ as a function of $\kappa$ behaves in this case like
 \[  \frac{1}{2}\kappa c^2 d + C\kappa^{-\frac{p-1}{1+2/d-p}}.  \]
Consequently, if we optimize the resulting bound for $\ov\la(p)$ over $c$, we get
\[ \ov\la(p)\leq \inf_{c\geq0} \left( \frac{1}{2}\kappa c d + 
C c^{-1}\kappa^{-\frac{p-1}{1+2/d-p}} \right)  = C^\prime 
\kappa^{\frac{1+1/d-p}{1+2/d-p}},  \]
 which implies the upper bound in \eqref{kappato0-2}. 
 
 In order to establish the lower bounds in \eqref{p-upperbound-2} and  \eqref{kappato0-2}, it suffices by the same reason as in (1) to take $b=0$. In this case, for fixed $\eps$ and $\kappa$, we bound \eqref{htildeintegral} from below by 
 \[
 C\frac{\tilde\al^{-2({1-(p-1)\frac{d}{2}})}}{1+\frac{2}{d}-p},
 \]
 where $C>0$ does not depend on $p$. As a result,
 \[
 \un\la(p)\geq 
 \left(\frac{C}{1+\frac{2}{d}-p}\right)^{\frac{1}{2(1-(p-1)d/2)}}, 
 \]
 which is the lower bound in \eqref{p-upperbound-2}. 
 For $\kappa\to0$, we repeat the argument given in the proof of Theorem~\ref{intfronts}, but use the truncation $\bone_{\{ g(1;t,x)>\eps \}}$ instead of $\bone_{\{ g(t,x)>\eps \}}$ in \eqref{mu-choice}. Hence, instead of $\tilde h$ in \eqref{h-tilde}, the function of interest is
 \[ h'(t)=\int_{|x|\geq \tilde\al t} g^p(t,x)\bone_{\{g(1;t,x)>\eps\}} \,\dd x. \]
 If we redo the calculations from \eqref{eq:htilde-lower1} to \eqref{htildeintegral}, then instead of \eqref{eq:def_R}, we should consider $R'=\tilde\al^2/(\kappa \eps^{2/d})$ so that in the end, we obtain exactly the same lower bound for $\int_0^\infty  h'(t)\,\dd t$ as in \eqref{htildeintegral}, but under the new condition $\tilde\al^2\eps^{-2/d}\geq 2\pi\kappa$.
Hence, we can make $\int_0^\infty  h'(t)\,\dd t$ arbitrarily large
if we take
\[\tilde \al = C\kappa^{\frac{1+1/d-p}{1+2/d-p}},\]
 and a large value for $C>0$. This choice of
$\tilde\al$ satisfies $\tilde\al^2\eps^{-2/d}\geq 2\pi\kappa$ for all 
$\kappa$ small enough, so the lower bound in \eqref{kappato0-2} follows. Note that at this 
part it is enough if $f(x) = O(e^{-c |x|})$ holds for some fixed $c >0$.
\eenu
\epr

\section*{Acknowledgements}

We would like to thank an anonymous referee for constructive comments, which in particular led to a more general statement in Theorem~\ref{compprin}.
CC acknowledges financial support from the Deutsche Forschungsgemeinschaft (project number KL 1041/7-1).
This research was initiated while PK held an Alexander von Humboldt postdoctoral 
fellowship at the Technical University of Munich. PK's research was further supported by the János Bolyai Research Scholarship of the Hungarian Academy of Sciences, and by the NKFIH grant FK124141.

\bibliographystyle{plain}
\bibliography{heat}

\begin{thebibliography}{10}

\bibitem{Ahn92}
H.S. Ahn, R.A. Carmona, and S.A. Molchanov.
\newblock Nonstationary {A}nderson model with {L}{\'e}vy potential.
\newblock In B.L. Rozovskii and R.B. Sowers, editors, {\em Stochastic Partial
  Differential Equations and Their Applications}, pages 1--11. Springer,
  Berlin, 1992.

\bibitem{Asmussen03}
S.~Asmussen.
\newblock {\em Applied Probability and Queues}.
\newblock Springer, New York, 2nd edition, 2003.

\bibitem{Balan16}
R.M. Balan and C.B. Ndongo.
\newblock Intermittency for the wave equation with {L}{\'e}vy white noise.
\newblock {\em Stat. Probab. Lett.}, 109:214--223, 2016.

\bibitem{Bertini95}
L.~Bertini and N.~Cancrini.
\newblock The stochastic heat equation: {Feynman}-{K}ac formula and
  intermittence.
\newblock {\em J. Stat. Phys.}, 78(5/6):1377--1401, 1995.

\bibitem{Besala63}
P.~Besala.
\newblock On solutions of {F}ourier's first problem for a system of non-linear
  parabolic equations in an unbounded domain.
\newblock {\em Ann. Polon. Math.}, 13(3):247--265, 1963.

\bibitem{Bichteler83}
K.~Bichteler and J.~Jacod.
\newblock Random measures and stochastic integration.
\newblock In G.~Kallianpur, editor, {\em Theory and Application of Random
  Fields}, pages 1--18. Springer, Berlin, 1983.

\bibitem{Carmona94}
R.A. Carmona and S.A. Molchanov.
\newblock {\em Parabolic Anderson Model and Intermittency}.
\newblock American Mathematical Society, Providence, RI, 1994.

\bibitem{Chen16}
B.~Chen, C.~Chong, and C.~Kl{\"u}ppelberg.
\newblock Simulation of stochastic {V}olterra equations driven by space--time
  {L}{\'e}vy noise.
\newblock In M.~Podolskij, R.~Stelzer, S.~Thorbj{\o}rnsen, and A.E.D. Veraart,
  editors, {\em The Fascination of Probability, Statistics and their
  Applications}, pages 209--229. Springer, Cham, Switzerland, 2016.

\bibitem{Chen15}
L.~Chen and R.C. Dalang.
\newblock Moments and growth indices for the nonlinear stochastic heat equation
  with rough initial conditions.
\newblock {\em Ann. Probab.}, 43(6):3006--3051, 2015.

\bibitem{Chong}
C.~Chong.
\newblock L\'evy-driven {V}olterra equations in space and time.
\newblock {\em J. Theor. Probab.}, 30(3):1014--1058, 2017.

\bibitem{Chong1}
C.~Chong.
\newblock Stochastic {PDE}s with heavy-tailed noise.
\newblock {\em Stoch. Process. Appl.}, 127(7):2262--2280, 2017.

\bibitem{Conus12}
D.~Conus and D.~Khoshnevisan.
\newblock On the existence and position of the farthest peaks of a family of
  stochastic heat and wave equations.
\newblock {\em Probab. Theory Relat. Fields}, 152(3):681--701, 2012.

\bibitem{Cranston05}
M.~Cranston, T.S. Mountford, and T.~Shiga.
\newblock Lyapunov exponent for the parabolic {A}nderson model with {L}{\'e}vy
  noise.
\newblock {\em Probab. Theory Relat. Fields}, 132(3):321--355, 2005.

\bibitem{Dalang09}
R.C. Dalang and C.~Mueller.
\newblock Intermittency properties in a hyperbolic {A}nderson model.
\newblock {\em Ann. Inst. Henri Poincar{\'e} Probab. Stat.}, 45(4):1150--1164,
  2009.

\bibitem{Dellacherie82}
C.~Dellacherie and P.-A. Meyer.
\newblock {\em Probabilities and Potential B}.
\newblock North-Holland, Amsterdam, 1982.

\bibitem{Doob53}
J.L. Doob.
\newblock {\em Stochastic Processes}.
\newblock Wiley, New York, 1953.

\bibitem{Foondun09}
M.~Foondun and D.~Khoshnevisan.
\newblock Intermittence and nonlinear stochastic partial differential
  equations.
\newblock {\em Electron. J. Probab.}, 14:548--568, 2009.

\bibitem{Foondun10}
M.~Foondun and D.~Khoshnevisan.
\newblock On the global maximum of the solution to a stochastic heat equation
  with compact-support initial data.
\newblock {\em Ann. Inst. Henri Poincar{\'e} Probab. Stat.}, 46(4):895--907,
  2010.

\bibitem{Hu16}
Y.~Hu, J.~Huang, and D.~Nualart.
\newblock On the intermittency front of stochastic heat equation driven by
  colored noises.
\newblock {\em Electron. Commun. Probab.}, 21(21):1--13, 2016.

\bibitem{Khos}
D.~Khoshnevisan.
\newblock {\em Analysis of Stochastic Partial Differential Equations}.
\newblock American Mathematical Society, Providence, RI, 2014.

\bibitem{MR}
C.~Marinelli and M.~R{\"o}ckner.
\newblock On maximal inequalities for purely discontinuous martingales in
  infinite dimensions.
\newblock In C.~Donati-Martin, A.~Lejay, and A.~Rouault, editors, {\em
  S\'eminaire de {P}robabilit\'es {XLVI}}, pages 293--315. Springer, Cham,
  Switzerland, 2014.

\bibitem{Mueller91}
C.~Mueller.
\newblock On the support of solutions to the heat equation with noise.
\newblock {\em Stoch. Stoch. Rep.}, 37(4):225--245, 1991.

\bibitem{SLB98}
E.~{Saint Loubert Bi{\'e}}.
\newblock {\'E}tude d'une {EDPS} conduite par un bruit poissonnien.
\newblock {\em Probab. Theory Relat. Fields}, 111(2):287--321, 1998.

\bibitem{Shiryaev96}
A.N. Shiryaev.
\newblock {\em Probability}.
\newblock Springer, New York, 2nd edition, 1996.

\bibitem{Veraar06}
M.C. Veraar.
\newblock {\em Stochastic Integration in Banach Spaces and Applications to
  Parabolic Evolution Equations}.
\newblock PhD thesis, Technical University of Delft, 2006.

\bibitem{Walsh86}
J.B. Walsh.
\newblock An introduction to stochastic partial differential equations.
\newblock In P.L. Hennequin, editor, {\em \'Ecole d'\'Et\'e de Probabilit\'es
  de Saint Flour XIV - 1984}, pages 265--439. Springer, Berlin, 1986.

\end{thebibliography}

\end{document}